\documentclass{amsart}
\usepackage{amsfonts,times,verbatim}
\usepackage{amssymb}
\usepackage{amsthm}
\usepackage{pstricks}

\numberwithin{equation}{section}
 
\newtheorem{Lemma}{Lemma}
\newtheorem{Prop}[Lemma]{Proposition}
\newtheorem{Crit}[Lemma]{Criterion}
\newtheorem{theorem}[Lemma]{Theorem}
\newtheorem{Cor}[Lemma]{Corollary}
\newtheorem{Conj}[Lemma]{Conjecture}

\theoremstyle{definition}
\newtheorem{Example}[Lemma]{Example}
\newtheorem{Def}[Lemma]{Definition}
\newtheorem{algorithm}[Lemma]{Algorithm}

\theoremstyle{remark}
\newtheorem{Rem}[Lemma]{Remark}

\numberwithin{Lemma}{section}

\newcommand{\bks}[1]{B^{#1,s}}
\newcommand{\cl}{\mathrm{cl}}
\newcommand{\cw}{\mathrm{cw}}
\newcommand{\dual}{*}
\newcommand{\es}{\varnothing}
\newcommand{\etil}[1]{\tilde{e}_{#1}}
\newcommand{\ftil}[1]{\tilde{f}_{#1}}
\newcommand{\hsig}{\check{\sigma}}
\newcommand{\inner}[2]{\langle #1\,,\,#2\rangle}
\newcommand{\La}{\Lambda}
\newcommand{\Le}{\hbox{\rotatedown{$\Gamma$}}}
\newcommand{\Lee}{\hbox{\rotatedown{L}}}
\newcommand{\Lt}{\om_{2}}
\newcommand{\la}{\lambda}
\newcommand{\lan}{\langle}
\newcommand\lev{\mbox{\sl lev}\,}
\newcommand{\ran}{\rangle}
\newcommand{\Lz}{\Lambda_{0}}
\newcommand{\Lo}{\Lambda_{1}}
\newcommand{\om}{\varpi}
\newcommand{\Path}{\mathcal{P}}
\newcommand{\stack}[2]{\genfrac{}{}{0pt}{}{#1}{#2}}
\newcommand{\Sup}{\mathcal{BC}}
\newcommand{\Tred}{T^{\#}}

\newcommand{\uq}[1]{U_{q}(\mathfrak{#1})}
\newcommand{\uqp}[1]{U_{q}'(\mathfrak{#1})}

\newcommand{\qbin}[2]{\genfrac{[}{]}{0pt}{}{#1}{#2}}
\newcommand{\qq}{\mathbb{Q}(q)}
\newcommand{\vare}{\varepsilon}
\newcommand{\wt}{\mathrm{wt}}
\newcommand{\T}{\mathcal{M}}
\newcommand{\tbt}{\tilde{B}^{2,2}}
\newcommand{\tbs}{\tilde{B}^{2,s}}
\newcommand{\iotil}{\check{\iota}}
\newcommand{\bcdual}{*_{\mathcal{BC}}}
\newcommand{\Z}{\mathbb{Z}}
\newcommand{\Mmin}[1]{\mathcal{M}_{\min}(#1)}

\begin{document}

\title{Finite-Dimensional Crystals $B^{2,s}$ for Quantum Affine Algebras 
of type $D_{n}^{(1)}$}

\author[A.~Schilling]{Anne Schilling}
\address{Department of Mathematics, University of California, One Shields
Avenue, Davis, CA 95616-8633, U.S.A.}
\email{anne@math.ucdavis.edu}
\urladdr{http://www.math.ucdavis.edu/$\sim$anne}

\author[P.~Sternberg]{Philip Sternberg}
\address{Department of Mathematics, University of California, One Shields
Avenue, Davis, CA 95616-8633, U.S.A.}
\email{sternberg@math.ucdavis.edu}
\urladdr{http://www.math.ucdavis.edu/$\sim$sternberg}

\thanks{\textit{Date:} August 2004; revised October 2005}
\thanks{Supported in part by the NSF grants DMS-0135345 and DMS-0200774.}


\begin{abstract} 
The Kirillov--Reshetikhin modules $W^{r,s}$ are finite-dimensional representations
of quantum affine algebras $U'_q(\mathfrak{g})$, labeled by a Dynkin node $r$ of
the affine Kac--Moody algebra $\mathfrak{g}$ and a positive integer $s$. 
In this paper we study the combinatorial structure of the crystal basis $B^{2,s}$ 
corresponding to $W^{2,s}$ for the algebra of type $D_n^{(1)}$. 
\end{abstract}

\maketitle

\section{Introduction}
Quantum algebras were introduced independently by Drinfeld
\cite{D:1985} and Jimbo \cite{J:1985} in their study of two
dimensional solvable lattice models in statistical mechanics.  Since
then quantum algebras have surfaced in many areas of mathematics and
mathematical physics, such as the theory of knots and links,
representation theory, and topological quantum field theory.  Of
special interest, in particular for lattice models and representation
theory, are finite-dimensional representations of quantum affine
algebras.  The irreducible finite-dimensional $U'_q(\mathfrak{g})$-modules 
for an affine Kac--Moody algebra $\mathfrak{g}$ were classified by Chari and 
Pressley~\cite{CP:1995,CP:1998} in terms of Drinfeld polynomials.  The
Kirillov--Reshetikhin modules $W^{r,s}$, labeled by a Dynkin node $r$
and a positive integer $s$, form a special class of these
finite-dimensional modules.  They naturally correspond to the weight
$s\om_r$, where $\om_r$ is the $r$-th fundamental weight of the underlying
finite algebra $\overline{\mathfrak{g}}$.

Kashiwara~\cite{Kash:1990, Kash:1991} showed that in the limit $q\to 0$ the highest-weight representations 
of the quantum algebra $U_q(\mathfrak{g})$ have very special bases, called
crystal bases. This construction makes it possible to study modules over quantum
algebras in terms of crystals graphs, which are purely combinatorial objects.
However, in general it is not yet known which finite-dimensional representations 
of affine quantum algebras have crystal bases and what their combinatorial structure
is. Recently, Hatayama et al.~\cite{HKOTT:2001,HKOTY:1999}
conjectured that crystal bases $B^{r,s}$ for the Kirillov--Reshetikhin modules 
$W^{r,s}$ exist. For type $A_{n}^{(1)}$, the crystals $B^{r,s}$ are known to 
exist~\cite{KKMMNN:1992}, and the explicit combinatorial crystal structure is also 
well-understood~\cite{Sh:2002}. Assuming that the crystals $B^{r,s}$ exist, their
structure for non-simply laced algebras can be described in terms of 
virtual crystals introduced in~\cite{OSS:2003,OSS:2003a}. The virtual crystal
construction is based on the following well-known algebra embeddings of non-simply 
laced into simply laced types:
\begin{eqnarray*}
    C_{n}^{(1)},A_{2n}^{(2)},A_{2n}^{(2)\dagger},D_{n+1}^{(2)} &
    \hookrightarrow & A_{2n-1}^{(1)} \\
    A_{2n-1}^{(2)},B_{n}^{(1)} & \hookrightarrow & D_{n+1}^{(1)} \\
    E_{6}^{(2)},F_{4}^{(1)} & \hookrightarrow & E_{6}^{(1)} \\
    D_{4}^{(3)},G_{2}^{(1)} & \hookrightarrow & D_{4}^{(1)} . 
\end{eqnarray*}

The main open problems in the theory of finite-dimensional affine crystals
are therefore the proof of the existence of $B^{r,s}$ and the combinatorial
structure of these crystals for types $D_{n}^{(1)}\;(n\geq 4)$ and
$E_{n}^{(1)}\;(n=6,7,8)$.  In this paper, we concentrate on type
$D_{n}^{(1)}$.  For irreducible representations corresponding to
multiples of the first fundamental weight (indexed by a one-row Young
diagram) or any single fundamental weight (indexed by a one-column
Young diagram) the crystals have been proven to exist and the structure is
known~\cite{KKMMNN:1992,Koga:1999}. In~\cite{HKOTT:2001,HKOTY:1999}, a conjecture 
is presented on the decomposition of $B^{r,s}$ as a crystal for the underlying 
finite algebra of type $D_n$. Specifically, as a type $D_n$ classical crystal
the crystals $B^{r,s}$ of type $D_n^{(1)}$ for $r\le n-2$ decompose as
\begin{equation*}
\bks{r}\cong \bigoplus_{\Lambda}B(\Lambda),
\end{equation*}
where the direct sum is taken over all weights $\Lambda$ for the finite algebra 
corresponding to  partitions obtained from 
an $r\times s$ rectangle by removing any number of $2\times 1$ vertical dominoes.  
Here $B(\Lambda)$ is the $U_q(D_n)$-crystal 
associated with the highest weight representation of highest weight $\Lambda$
(see~\cite{KN:1994}).  In the sequel, we consider the case $r=2$, for which the 
above direct sum specializes to
\begin{equation}\label{eq:decomp}
\bks{2}\cong \bigoplus_{k=0}^{s}B(k\Lt),
\end{equation}
where once again the summands in the right hand side of the equation are crystals for the finite algebra.
Our approach to study the combinatorics of $B^{2,s}$ is as follows.
First, we introduce tableaux of shape $(s,s)$ to define a
$U_q(D_n)$-crystal whose vertices are in bijection with the classical
tableaux from the direct sum decomposition~\eqref{eq:decomp}.  Using
the automorphism of the $D_n^{(1)}$ Dynkin diagram which interchanges
nodes 0 and 1, we define the unique action of $\ftil{0}$ and
$\etil{0}$ which makes this crystal into a perfect crystal $\tbs$ of
level $s$ with an energy function.  (See sections \ref{subsec:perfect} and \ref{subsec:energy} for definitions of these terms.)

Assuming the existence of the crystal $B^{r,s}$, the main result of our paper
states that our combinatorially constructed crystal $\tbs$ is the unique perfect
crystal of level $s$ with the classical decomposition~\eqref{eq:decomp} with 
a given energy function. More precisely:
\begin{theorem} \label{thm:main} 
If $B^{2,s}$ exists with the properties as in Conjecture~\ref{conj:Brs}, 
then $\tbs\cong B^{2,s}$.
\end{theorem}

This is the first step in confirming Conjecture 2.1 of \cite{HKOTT:2001},
which states that as modules over the embedded classical quantum
group, $W^{2,s}$ decomposes as $\bigoplus_{k=0}^{s}V(k\Lt)$, where
$V(\Lambda)$ is the classical module with highest weight $\Lambda$,
$W^{2,s}$ has a crystal basis, and this crystal is a perfect crystal of level
$s$.

The paper is structured as follows. In section~\ref{sec:review} the definition
of quantum algebras, crystal bases and perfect crystals is reviewed.
Section~\ref{sec:typeD} is devoted to crystals and the plactic monoid
of type $D_n$. The properties of $B^{2,s}$ of type $D_n^{(1)}$ as conjectured 
in~\cite{HKOTT:2001} are given in Conjecture~\ref{conj:Brs}. In section~\ref{sec:decomp} 
the set underlying $\tbs$ is constructed
in terms of tableaux of shape $(s,s)$ obeying certain conditions. It is shown
that this set is in bijection with the union of sets appearing on the right
hand side of~\eqref{eq:decomp}. The branching component graph is introduced 
in section~\ref{sec:bc}, which is used in section~\ref{sec:aff} to define
$\etil{0}$ and $\ftil{0}$ on $\tbs$.  
This makes $\tbs$ into an affine crystal.
It is shown in section~\ref{sec:perfect} that $\tbs$ is perfect and that $\tbs$ is 
the unique perfect crystal having the classical decomposition~\eqref{eq:decomp}
with the appropriate energy function. 
This proves in particular Theorem~\ref{thm:main}. Finally, we end in 
section~\ref{sec:discussion} with some open problems.

\subsection*{Acknowledgments}
We would like to thank the Max-Planck-Institut f{\"u}r Mathematik in Bonn
and the Research Institute for Mathematical Sciences in Kyoto where part of 
this work was carried out. 
We would also like to thank Mark Haiman, Tomoki Nakanishi, Masato Okado and
Mark Shimozono for helpful conversations.
AS would like to thank the University of British Columbia for their hospitality 
during her forced exile from the US.

\section{Review of quantum groups and crystal bases}\label{sec:review}

\subsection{Quantum groups}
For $n\in\mathbb{Z}$ and a formal parameter $q$, we use the notation
\begin{displaymath}
    [n]_{q}=\frac{q^{n}-q^{-n}}{q-q^{-1}},\quad
    [n]_{q}!=\prod_{k=1}^{n}[k]_{q},\textrm{ and }
    \qbin{m}{n}_q=\frac{[m]_{q}!}{[n]_{q}![m-n]_{q}!}.
\end{displaymath}
These are all elements of $\qq$, called the $q$-integers, 
$q$-factorials, and $q$-binomial coefficients, respectively.

Let $\mathfrak{g}$ be an arbitrary Kac-Moody Lie algebra with Cartan datum
$(A,\Pi,\Pi^{\vee},P,P^{\vee})$ and a Dynkin diagram indexed by $I$.
Here $A=(a_{ij})_{i,j\in I}$ is the Cartan matrix, $P$ and $P^\vee$
are the weight lattice and dual weight lattice, respectively,
$\Pi=\{\alpha_i \mid i\in I\}$ is the set of simple roots and
$\Pi^\vee=\{h_i \mid i\in I\}$ is the set of simple coroots.
Furthermore, let $\{s_i \mid i\in I\}$ be the entries of the diagonal 
symmetrizing matrix of $A$ and define $q_{i}=q^{s_{i}}$ and $K_{i}=q^{s_{i}h_{i}}$.
Then the quantum enveloping algebra $\uq{g}$ is the associative $\qq$-algebra 
generated by $e_{i}$ and $f_{i}$ for $i\in I$, and $q^{h}$ for $h\in P^{\vee}$, 
with the following relations (see e.g. \cite[Def. 3.1.1]{HK:2002}):
\begin{enumerate}
    \item $q^{0}=1$, $q^{h}q^{h'}=q^{h+h'}$ for all $h,h'\in P^{\vee}$,
    \item $q^{h}e_{i}q^{-h}=q^{\alpha_{i}(h)}e_{i}$ for all $h\in P^{\vee}$,
    \item $q^{h}f_{i}q^{-h}=q^{\alpha_{i}(h)}f_{i}$ for all $h\in P^{\vee}$,
    \item $e_{i}f_{j}-f_{j}e_{i}=\delta_{ij}\frac{K_{i}-K_{i}^{-1}}{q_{i}-q_{i}^{-1}}$ 
          for $i,j\in I$,
    \item $\sum_{k=0}^{1-a_{ij}}(-1)^{k}\qbin{1-a_{ij}}{k}_{q_{i}}e_{i}^{1-a_{ij}-k}
          e_{j}e_{i}^{k}=0$ for all $i\neq j$,
    \item $\sum_{k=0}^{1-a_{ij}}(-1)^{k}\qbin{1-a_{ij}}{k}_{q_{i}}f_{i}^{1-a_{ij}-k}
          f_{j}f_{i}^{k}=0$ for all $i\neq j$.
\end{enumerate}

\subsection{Crystal bases}\label{sec:crystal bases}
The quantum algebra $\uq{g}$ can be viewed as a $q$-deformation of the
universal enveloping algebra $U(\mathfrak{g})$ of $\mathfrak{g}$.
Lusztig~\cite{Lus:1988} showed that the integrable highest weight representations
of $U(\mathfrak{g})$ can be deformed to $U_q(\mathfrak{g})$ representations
in such a way that the dimension of the weight spaces are invariant under
the deformation, provided $q\neq0$ and $q^{k}\neq1$ for all $k\in\mathbb{Z}$
(see also~\cite{HK:2002}). Let $M$ be a $U_q(\mathfrak{g})$-module and
$R$ the subset of all elements in $\qq$ which are regular at $q=0$.
Kashiwara~\cite{Kash:1990,Kash:1991} introduced Kashiwara operators
$\etil{i}$ and $\ftil{i}$ as certain linear combinations of powers of
$e_i$ and $f_i$. A crystal lattice $\mathcal{L}$ is a free $R$-submodule
of $M$ that generates $M$ over $\qq$, has the same weight decomposition
and has the property that $\etil{i}\mathcal{L}\subset \mathcal{L}$ and
$\ftil{i}\mathcal{L}\subset \mathcal{L}$ for all $i\in I$. The passage
from $\mathcal{L}$ to the quotient $\mathcal{L}/q\mathcal{L}$ is referred
to as taking the crystal limit. A crystal basis is a $\mathbb{Q}$-basis
of $\mathcal{L}/q\mathcal{L}$ with certain properties.

Axiomatically, we may define a $U_q(\mathfrak{g})$-crystal as a nonempty set $B$
equipped with maps $\wt:B\rightarrow P$ and
$\etil{i},\ftil{i}:B\rightarrow B\cup\{\es\}$ for all $i\in I$,
satisfying
\begin{align}
\label{eq:e-f}
\ftil{i}(b)=b' &\Leftrightarrow \etil{i}(b')=b
\text{ if $b,b'\in B$} \\
\wt(\ftil{i}(b))&=\wt(b)-\alpha_i \text{ if $\ftil{i}(b)\in B$} \\
\label{eq:string length}
\inner{h_i}{\wt(b)}&=\varphi_i(b)-\epsilon_i(b).
\end{align}
Here for $b \in B$
\begin{equation*}
\begin{split}
\epsilon_i(b)&= \max\{n\ge0\mid \etil{i}^n(b)\not=\es \} \\
\varphi_i(b) &= \max\{n\ge0\mid \ftil{i}^n(b)\not=\es \}.
\end{split}
\end{equation*}
(It is assumed that $\varphi_i(b),\epsilon_i(b)<\infty$ for all
$i\in I$ and $b\in B$.) A $U_q(\mathfrak{g})$-crystal $B$ can be viewed
as a directed edge-colored graph (the crystal graph) whose
vertices are the elements of $B$, with a directed edge from $b$ to
$b'$ labeled $i\in I$, if and only if $\ftil{i}(b)=b'$.

Let $B_1$ and $B_2$ be $U_q(\mathfrak{g})$-crystals. The Cartesian product
$B_2\times B_1$ can also be endowed with the structure of a $U_q(\mathfrak{g})$-crystal.
The resulting crystal is denoted by $B_2\otimes B_1$ and its elements $(b_2,b_1)$ are 
written $b_2\otimes b_1$.
(The reader is warned that our convention is opposite to that of 
Kashiwara~\cite{Kash:1995}). For $i\in I$ and $b=b_2\otimes b_1\in B_2\otimes B_1$,
we have $\wt(b)=\wt(b_1)+\wt(b_2)$,
\begin{equation} \label{eq:f on two factors}
\ftil{i}(b_2\otimes b_1) = \begin{cases} \ftil{i}(b_2)\otimes b_1
& \text{if $\epsilon_i(b_2)\ge \varphi_i(b_1)$} \\
b_2\otimes \ftil{i}(b_1) & \text{if $\epsilon_i(b_2)<\varphi_i(b_1)$}
\end{cases}
\end{equation}
and
\begin{equation} \label{eq:e on two factors}
\etil{i}(b_2\otimes b_1) = \begin{cases} \etil{i}(b_2) \otimes b_1 &
\text{if $\epsilon_i(b_2)>\varphi_i(b_1)$} \\
b_2\otimes \etil{i}(b_1) & \text{if $\epsilon_i(b_2)\le \varphi_i(b_1)$.}
\end{cases}
\end{equation}
Combinatorially, this action of $\ftil{i}$ and $\etil{i}$ on tensor
products can be described by the signature rule.
The $i$-signature of $b$ is the word consisting of the symbols $+$ and $-$ 
given by
\begin{equation*}
\underset{\text{$\varphi_i(b_2)$ times}}{\underbrace{-\dotsm-}}
\quad \underset{\text{$\epsilon_i(b_2)$
times}}{\underbrace{+\dotsm+}} \quad 
\underset{\text{$\varphi_i(b_1)$ times}}{\underbrace{-\dotsm-}}
\quad \underset{\text{$\epsilon_i(b_1)$
times}}{\underbrace{+\dotsm+}} .
\end{equation*}
The reduced $i$-signature of $b$ is the subword of the
$i$-signature of $b$, given by the repeated removal of adjacent
symbols $+-$ (in that order); it has the form
\begin{equation*}
\underset{\text{$\varphi$ times}}{\underbrace{-\dotsm-}} \quad
\underset{\text{$\epsilon$ times}}{\underbrace{+\dotsm+}}.
\end{equation*}
If $\varphi=0$ then $\ftil{i}(b)=\es$; otherwise $\ftil{i}$ acts
on the tensor factor corresponding to the rightmost symbol $-$ in the 
reduced $i$-signature of $b$. Similarly, if $\epsilon=0$ then
$\etil{i}(b)=\es$; otherwise $\etil{i}$ acts on the leftmost symbol $+$ 
in the reduced $i$-signature of $b$. From this it is clear that
\begin{equation*}
\begin{split}
&\varphi_{i}(b_{2}\otimes 
b_{1})=\varphi_{i}(b_{2})+\mathrm{max}(0,\varphi_{i}(b_{1})-\vare_{i}(b_{2})) ,\\
&\vare_{i}(b_{2}\otimes
b_{1})=\vare_{i}(b_{1})+\mathrm{max}(0,-\varphi_{i}(b_{1})+\vare_{i}(b_{2})).
\end{split}
\end{equation*}

\subsection{Perfect crystals} \label{subsec:perfect}
Of particular interest is a class of crystals called perfect crystals, which are crystals for affine algebras satisfying a set of very special properties.  These properties ensure that perfect crystals can be used to construct the path realization of highest weight modules~\cite{KKMMNN:1992a}. To define them, we need a few preliminary definitions.

Recall that $P$ denotes the weight lattice of a Kac-Moody algebra $\mathfrak{g}$; for the remainder of this section, $\mathfrak{g}$ is of affine type.  The center of $\mathfrak{g}$ is one-dimensional and is generated by the 
canonical central element $c=\sum_{i\in I}a^\vee_i h_i$, where the $a^\vee_i$ are the numbers on the nodes of the Dynkin diagram of the algebra dual to $\mathfrak{g}$ given in Table Aff of~\cite[section 4.8]{Kac:1990}.  Moreover, the imaginary
roots of $\mathfrak{g}$ are nonzero integral multiples of the null root
$\delta=\sum_{i\in I} a_i \alpha_i$, where the $a_i$ are the numbers on the nodes of the Dynkin diagram of $\mathfrak{g}$ given in Table Aff of~\cite{Kac:1990}.  Define $P_{\cl}=P/\mathbb{Z}\delta$,
$P_{\cl}^{+}=\{\la\in P_{\cl}|\lan h_{i},\la\ran\geq 0\textrm{ for all 
}i\in I\}$, and $\uqp{g}$ to be the quantum enveloping algebra with the Cartan datum
$(A,\Pi,\Pi^{\vee},P_{\cl},P_{\cl}^{\vee})$.

Define the set of level $\ell$ weights to be
$(P_{\cl}^{+})_{\ell}=\{\la\in P_{\cl}^{+}|\lan c,\la \ran=\ell\}$.  For
a crystal basis element $b\in B$, define
\begin{equation*}
\vare(b)=\sum_{i\in I}\vare_{i}(b)\Lambda_{i} \quad \text{and} \quad
\varphi(b)=\sum_{i\in I}\varphi_{i}(b)\Lambda_{i},
\end{equation*}
where $\Lambda_i$ is the $i$-th fundamental weight of $\mathfrak{g}$.
Finally, for a crystal basis $B$, we define $B_{\min}$ to be the set of crystal basis
elements $b$ such that $\lan c,\vare(b)\ran$ is minimal over $b\in B$.

\begin{Def} \label{def:perfect}
A crystal $B$ is a perfect crystal of level $\ell$ if:
\begin{enumerate}
    \item $B\otimes B$ is connected;
    \item there exists $\la\in P_{\cl}$ such that $\wt(B)\subset
    \la+\sum_{i\neq0}\mathbb{Z}_{\leq0}\alpha_{i}$ and
    $\#(B_{\la})=1$;
    \item there is a finite-dimensional irreducible $\uqp{g}$-module $V$ with a
    crystal base whose crystal graph is isomorphic to $B$;
    \item for any $b\in B$, we have $\langle c,\vare(b)\rangle \geq
    \ell$;
    \item the maps $\vare$ and $\varphi$ from $B_{\min}$ to
    $(P_{\cl}^{+})_{\ell}$ are bijective.
\end{enumerate}
We use the notation $\lev(B)$ to indicate the level of the perfect 
crystal $B$.
\end{Def}

\subsection{Energy function} \label{subsec:energy}

The existence of an affine crystal structure usually provides an energy function.
Let $B_1$ and $B_2$ be finite $U'_q(\mathfrak{g})$-crystals. Then
following~\cite[Section 4]{KKMMNN:1992a} we have:
\begin{enumerate}
\item There is a unique isomorphism of $U'_q(\mathfrak{g})$-crystals
$R=R_{B_2,B_1}:B_2\otimes B_1\rightarrow B_1\otimes B_2$.
\item There is a function $H=H_{B_2,B_1}:B_2\otimes
B_1\rightarrow\Z$, unique up to global additive constant, such
that $H$ is constant on classical components and, for all $b_2\in B_2$
and $b_1\in B_1$, if $R(b_2\otimes b_1)=b_1'\otimes b_2'$, then
\begin{equation} \label{eq:local energy}
  H(\etil{0}(b_2\otimes b_1))=
  H(b_2\otimes b_1)+
  \begin{cases}
    -1 & \text{if $\vare_0(b_2)>\varphi_0(b_1)$ and
    $\vare_0(b_1')>\varphi_0(b_2')$} \\
    1 & \text{if $\vare_0(b_2)\le\varphi_0(b_1)$ and
    $\vare_0(b_1')\le \varphi_0(b_2')$} \\
    0 & \text{otherwise.}
  \end{cases}
\end{equation}
\end{enumerate}

We shall call the maps $R$ and $H$ the local isomorphism and
local energy function on $B_2\otimes B_1$. The pair $(R,H)$ is
called the combinatorial $R$-matrix.

Let $u(B_1)$ and $u(B_2)$ be extremal vectors of $B_1$ and $B_2$, respectively
(see \cite{Kash:2002} for a definition of extremal vectors). Then
\begin{equation*}
  R(u(B_2)\otimes u(B_1)) = u(B_1)\otimes u(B_2).
\end{equation*}
It is convenient to normalize the local energy function $H$ by
requiring that
\begin{equation*}
  H(u(B_2)\otimes u(B_1)) = 0.
\end{equation*}
With this convention it follows by definition that
\begin{equation*}
H_{B_1,B_2} \circ R_{B_2,B_1} = H_{B_2,B_1}
\end{equation*}
as operators on $B_2\otimes B_1$.

We wish to define an energy function
$D_B:B\rightarrow \Z$ for tensor products of perfect crystals
of the form $B^{r,s}$~\cite[Section 3.3]{HKOTT:2001}. 
Let $B=B^{r,s}$ be perfect. Then there exists
a unique element $b^\natural\in B$ such that $\varphi(b^\natural)=\lev(B)\La_0$.
Define $D_B:B\rightarrow\Z$ by
\begin{equation}\label{eq:D single}
  D_B(b) = H_{B,B}(b\otimes b^\natural) - H_{B,B}(u(B)\otimes b^\natural).
\end{equation}

The intrinsic energy $D_B$ for the $L$-fold tensor product
$B=B_L\otimes\dotsm\otimes B_1$ where $B_j=B^{r_j,s_j}$ is given by
\begin{equation*}
D_B = \sum_{1\le i<j\le L} H_i R_{i+1}R_{i+2}\dotsm R_{j-1}
 + \sum_{j=1}^L D_{B_j} R_1R_2\dotsm R_{j-1},
\end{equation*}
where $H_i$ and $R_i$ are the local energy function and 
$R$-matrix on the $i$-th and $i+1$-th tensor factor, respectively.

\section{Crystals and plactic monoid of type $D$}\label{sec:typeD}
{}From now on we restrict our attention to the finite Lie algebra of
type $D_n$ and the affine Kac-Moody algebra of type $D_n^{(1)}$.
Denote by $I=\{0,1,\ldots,n\}$ the index set of the Dynkin diagram for
$D_n^{(1)}$ and by $J=\{1,2,\ldots,n\}$ the Dynkin diagram for type
$D_n$.

\subsection{Dynkin data}
For type $D_n$, the simple roots are
\begin{equation}\label{eq:alpha}
\begin{split}
\alpha_i&=\varepsilon_i-\varepsilon_{i+1} \qquad \text{for $1\le i<n$}\\
\alpha_n&=\varepsilon_{n-1}+\varepsilon_n
\end{split}
\end{equation}
and the fundamental weights are
\begin{equation*}
\begin{aligned}
\om_i&=\varepsilon_1+\cdots+\varepsilon_i &\text{for $1\le i\le n-2$}\\
\om_{n-1}&=(\varepsilon_1+\cdots+\varepsilon_{n-1}-\varepsilon_n)/2 &\\
\om_n&=(\varepsilon_1+\cdots+\varepsilon_{n-1}+\varepsilon_n)/2&
\end{aligned}
\end{equation*}
where $\varepsilon_i\in \mathbb{Z}^n$ is the $i$-th unit standard vector.
The central element for $D_n^{(1)}$ is given by
$$c=h_{0}+h_{1}+2h_{2}+\cdots+2h_{n-2}+h_{n-1}+h_{n}.$$

\subsection{Classical crystals}\label{sec:class crys}
Kashiwara and Nakashima~\cite{KN:1994} described the crystal structure of all 
classical highest weight crystals $B(\Lambda)$ of highest weight $\Lambda$ explicitly. 
For the special case $B(k\Lt)$ as occuring in \eqref{eq:decomp} each crystal element
can be represented by a tableau of shape $\lambda=(k,k)$ on the partially ordered 
alphabet
\begin{displaymath}
1<2<\cdots < n-1 < \stack{n}{\bar{n}} < \overline{n-1}<\cdots 
\bar{2}<\bar{1}
\end{displaymath}
such that the following conditions hold~\cite[page 202]{HK:2002}:
\begin{Crit} \label{legal}
{}~
\begin{enumerate}
    \item \label{legal:1} If $ab$ is in the filling, then $a\leq b$;
    \item \label{legal:2} If $\stack{a}{b}$ is in the filling, then $b\nleq a$;
    \item \label{legal:3} No configuration of the form $\stack{a}{\phantom{\bar{a}}}
      \stack{a}{\bar{a}}$ or $\stack{a}{\bar{a}}\stack{\phantom{a}}{\bar{a}}$ appears;
    \item \label{legal:4} No configuration of the form 
      $\stack{n-1}{n}\ldots\stack{n}{\overline{n-1}}$ or 
      $\stack{n-1}{\bar{n}}\ldots\stack{\bar{n}}{\overline{n-1}}$ appears;
    \item \label{legal:5} No configuration of the form $\stack{1}{\bar{1}}$ appears.
\end{enumerate}
\end{Crit}
Note that for $k\geq 2$, condition $5$ follows from conditions $1$ and 
$3$.

Also, observe that the conditions given in~\cite{HK:2002} apply only
to adjacent columns, not to non-adjacent columns as in condition 4
above.  However, Criterion \ref{legal} is unchanged by replacing
condition 4 with the following:
\begin{enumerate}
    \item[(4a)] {\it No configuration of the form }
    $\stack{n-1}{n}\stack{n}{\overline{n-1}}$ {\it or}
    $\stack{n-1}{\bar{n}}\stack{\bar{n}}{\overline{n-1}}$ {\it
    appears.}
\end{enumerate}
    
To see this equivalence, observe that by conditions $1$ and $2$ the
only columns that can appear between $\stack{n-1}{n}$ and
$\stack{n}{\overline{n-1}}$ are $\stack{n-1}{n}$,
$\stack{n-1}{\overline{n-1}}$, and $\stack{n}{\overline{n-1}}$, and
they must appear in that order from left to right.  If a column of the
form $\stack{n-1}{\overline{n-1}}$ appears, we have a configuration of
the form $\stack{n-1}{\phantom{\overline{n-1}}}
\stack{n-1}{\overline{n-1}}$, which is forbidden by condition $3$.  On
the other hand, if no column of the form $\stack{n-1}{\overline{n-1}}$
appears, the columns $\stack{n-1}{n}$ and $\stack{n}{\overline{n-1}}$
are adjacent, which is disallowed by condition 4a.

The crystal $B(\om_1)$ is described pictorially by the crystal graph:
\vspace{.5in}
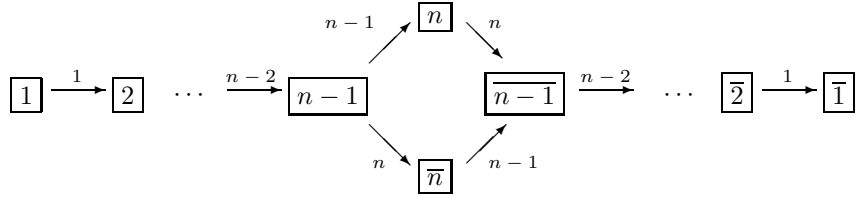
\begin{center}\setlength{\unitlength}{1.8cm}
    \begin{picture}(0,0)
	\put(-3.2,0){\framebox{$1$}}
	\put(-2.9,.1){\vector(1,0){.4}}
	\put(-2.75,.16){\scriptsize{$1$}}
	\put(-2.45,0){\framebox{$2$}}
	\put(-2,.02){$\cdots$}
	\put(-1.6,.1){\vector(1,0){.4}}
	\put(-1.61,.16){\scriptsize{$n-2$}}
	\put(-1.15,0){\framebox{$n-1$}}
	\put(-.55,.3){\vector(1,1){.3}}
	\put(-.88,.56){\scriptsize{$n-1$}}
	\put(-.55,-.15){\vector(1,-1){.3}}
	\put(-.53,-.47){\scriptsize{$n$}}
	\put(-.19,.6){\framebox{$n$}}
	\put(.17,.6){\vector(1,-1){.26}}
	\put(.33,.56){\scriptsize{$n$}}
	\put(-.19,-.6){\framebox{$\overline{n}$}}
	\put(.17,-.44){\vector(1,1){.29}}
	\put(.33,-.47){\scriptsize{$n-1$}}
	\put(0.3,0){\framebox{$\overline{n-1}$}}
	\put(1.,.1){\vector(1,0){.4}}
	\put(1.01,.16){\scriptsize{$n-2$}}
	\put(1.62,.02){$\cdots$}
	\put(2.05,0){\framebox{$\overline{2}$}}
	\put(2.35,.1){\vector(1,0){.4}}
	\put(2.50,.16){\scriptsize{$1$}}
	\put(2.8,0){\framebox{$\overline{1}$}}	
    \end{picture}
\end{center}
\vspace{.5in} 
For a tableau $T=\begin{array}{c} a_{1}\\ b_{1} \end{array} 
\cdots \begin{array}{c} a_{k}\\b_{k} \end{array} \in B(k\Lt)$,
the action of the Kashiwara operators $\ftil{i}$ and $\etil{i}$ 
is defined as follows.
Consider the column word $w_{T}=b_{1}a_{1}\cdots b_{k}a_{k}$ and
view this word as an element in $B(\om_1)^{\otimes 2k}$. Then
$\ftil{i}$ and $\etil{i}$ act by the tensor product rule as defined in
section~\ref{sec:crystal bases}.

\begin{Example}
    Let $n=4$.  Then the tableau
    $$
    T=\left.
    \begin{array}{|c|c|c|c|c|}
	\hline
	1 & 2 & 4 & \bar{3} & \bar{3}  \\
	\hline
	3 & \bar{4} & \bar{4} & \bar{2} & \bar{1}  \\
	\hline 
    \end{array}
    \right.
    $$
    has column word
    $w_{T}=31\bar{4}2\bar{4}4\bar{2}\bar{3}\bar{1}\bar{3}$.  The
    $2$-signature of $T$ is $+-+--$, derived from the subword
    $32\bar{2}\bar{3}\bar{3}$, and the reduced $2$-signature is a single
    $-$.  Therefore,
    $$
    \ftil{2}(T)=\left.
    \begin{array}{|c|c|c|c|c|}
	\hline
	1 & 2 & 4 & \bar{3} & \bar{2}  \\
	\hline
	3 & \bar{4} & \bar{4} & \bar{2} & \bar{1}  \\
	\hline 
    \end{array}
    \right. ,
    $$
    since the rightmost $-$ in the reduced $2$-signature of $T$ comes from
    the northeastmost $\bar{3}$.  The $4$-signature of $T$ is $-++-++$,
    derived from the subword $3\bar{4}\bar{4}4\bar{3}\bar{3}$, and the
    reduced $4$-signature is $-+++$, from the subword
    $3\bar{4}\bar{3}\bar{3}$.  This tells us that
    $$
    \ftil{4}(T)=\left.
    \begin{array}{|c|c|c|c|c|}
	\hline
	1 & 2 & 4 & \bar{3} & \bar{3}  \\
	\hline
	\bar{4} & \bar{4} & \bar{4} & \bar{2} & \bar{1}  \\
	\hline 
    \end{array}
    \right.
    \quad\textrm{ and }\quad
    \etil{4}(T)=\left.
    \begin{array}{|c|c|c|c|c|}
	\hline
	1 & 2 & 4 & \bar{3} & \bar{3}  \\
	\hline
	3 & 3 & \bar{4} & \bar{2} & \bar{1}  \\
	\hline 
    \end{array}
    \right. .
    $$    
\end{Example}

\subsection{Dual crystals}\label{sec:dual crystal}
Let $\omega_0$ be the longest element in the Weyl group of $D_n$.
The action of $\omega_0$ on the weight lattice $P$ of $D_n$ is given by
\begin{equation*}
\begin{split}
\omega_0(\om_i) &= -\om_{\tau(i)}\\
\omega_0(\alpha_i) &= -\alpha_{\tau(i)}
\end{split}
\end{equation*}
where $\tau:J\to J$ is the identity if $n$ is even and interchanges
$n-1$ and $n$ and fixes all other Dynkin nodes if $n$ is odd. 

There is a unique involution $\dual:B\to B$, called the dual map,
satisfying
\begin{equation*}
\begin{split}
\wt(b^\dual) &= \omega_0 \wt(b)\\
\etil{i}(b)^\dual &= \ftil{\tau(i)}(b^\dual)\\
\ftil{i}(b)^\dual &= \etil{\tau(i)}(b^\dual).
\end{split}
\end{equation*}
The involution $\dual$ sends the highest weight vector $u\in B(\La)$ to the 
lowest weight vector (the unique vector in $B(\La)$ of weight $\omega_0(\La)$).
We have
\begin{equation*}
(B_1\otimes B_2)^\dual \cong B_2\otimes B_1
\end{equation*}
with $(b_1\otimes b_2)^\dual \mapsto b_2^\dual \otimes b_1^\dual$.

Explicitly, on $B(\om_1)$ the involution $\dual$ is given by 
\begin{equation*}
 i \longleftrightarrow \overline{i}
\end{equation*}
except for $i=n$ with $n$ odd in which case $n \leftrightarrow n$ and 
$\overline{n} \leftrightarrow \overline{n}$. 
For $T\in B(\La)$ the dual $T^{\dual}$ is obtained by applying
the $\dual$ map defined for $B(\om_1)$ to each of the letters of 
$w_{T}^{\textrm{rev}}$ (the reverse column word of $T$), and then rectifying 
the resulting word.

\begin{Example}
If
\begin{equation*}
T=\begin{array}{|c|c|c|} \hline 1&1&2\\ \hline \overline{3} & \multicolumn{2}{c}{}\\
\cline{1-1}\end{array}
\in B(2\om_1+\om_2)
\end{equation*}
we have
\begin{equation*}
T^{*}=\begin{array}{|c|c|c|} \hline 3&\overline{1}&\overline{1}\\
\hline \overline{2} & \multicolumn{2}{c}{}\\ \cline{1-1}\end{array}.
\end{equation*}
\end{Example}

\subsection{Plactic monoid of type $D$}\label{sec:plac}
The plactic monoid for type $D$ is the free monoid generated by
$\{1,\ldots,n,\bar{n},\ldots,\bar{1}\}$, modulo certain relations
introduced by Lecouvey~\cite{L}. Note that we write our words in the reverse
order compared to~\cite{L}. A column word $C=x_L x_{L-1} \cdots x_1$
is a word such that $x_{i+1}\not\le x_i$ for $i=1,\ldots,L-1$.  Note
that the letters $n$ and $\bar{n}$ are the only letters that may
appear more than once in $C$.  Let $z\le n$ be a letter in $C$.  Then
$N(z)$ denotes the number of letters $x$ in $C$ such that $x\le z$ or
$x\ge \bar{z}$.  A column $C$ is called admissible if $L\le n$ and
for any pair $(z,\bar{z})$ of letters in $C$ with $z\le n$ we have
$N(z)\le z$.  The Lecouvey $D$ equivalence relations are given by:
\begin{enumerate}
    \item  If $x\neq \bar{z}$, then
    \begin{displaymath}
	xzy\equiv zxy \textrm{ for }x\leq y<z
	\textrm{ and }
	yzx\equiv yxz \textrm{ for }x<y\leq z	.
    \end{displaymath}
    \item  If $1<x<n$ and $x\leq y\leq\bar{x}$, then 
    \begin{displaymath}
	(x-1)(\overline{x-1})y\equiv\bar{x}xy\textrm{ and }
	y\bar{x}x\equiv y(x-1)(\overline{x-1}).
    \end{displaymath}
    \item  If $x\leq n-1$, then
    \begin{displaymath}
	\left\{
	\begin{array}{c}
	    n\bar{x}\bar{n}\equiv n\bar{n}\bar{x}\\
	    \bar{n}\bar{x}n\equiv \bar{n}n\bar{x}\\
	\end{array}
	\right.
	\textrm{ and }
	\left\{
	\begin{array}{c}
	    xn\bar{n}\equiv nx\bar{n}\\
	    x\bar{n}n\equiv \bar{n}xn\\
	\end{array}
	\right. .
    \end{displaymath}
    \item 
    \begin{displaymath}
	\left\{
	\begin{array}{c}
	    \bar{n}\bar{n}n\equiv\bar{n}(n-1)(\overline{n-1})\\
	    nn\bar{n}\equiv n(n-1)(\overline{n-1})\\
	\end{array}
	\right.
	\textrm{ and }
	\left\{
	\begin{array}{c}
	    (n-1)(\overline{n-1})\bar{n}\equiv n\bar{n}\bar{n}\\
	    (n-1)(\overline{n-1})n\equiv \bar{n}nn\\
	\end{array}
	\right. .
    \end{displaymath}
    \item 
    Consider $w$ a non-admissible column word each strict factor of 
    which is admissible.  Let $z$ be the lowest unbarred letter such 
    that the pair $(z,\bar{z})$ occurs in $w$ and $N(z)>z$.  Then 
    $w\equiv\tilde{w}$ is the column word obtained by erasing the 
    pair $(z,\bar{z})$ in $w$ if $z<n$, by erasing a pair 
    $(n,\bar{n})$ of consecutive letters otherwise.
\end{enumerate}
This monoid gives us a bumping algorithm similar to the Schensted bumping
algorithm.  It is noted in~\cite{L} that a general type $D$ sliding algorithm, 
if one exists, would be very complicated.  However, for tableaux with no more than 
two rows, these relations provide us with the following relations on 
subtableaux:
\begin{enumerate}

    \item  If $x\neq \bar{z}$, then
    \begin{eqnarray*}
	\left.
	\begin{array}{c|c|c|}
	    \cline{3-3}
	     \multicolumn{1}{c}{} & & y \\
	    \cline{2-3}
	    & x & z  \\
	    \cline{2-3}
	\end{array}
	\equiv
	\begin{array}{c|c|c|}
	    \cline{2-3}
	      & x& y \\
	    \cline{2-3}
	    \multicolumn{1}{c}{}&  & z  \\
	    \cline{3-3}
	\end{array}	
	\equiv
	\begin{array}{|c|c|c}
	    \cline{1-2}
	     x &  y \\
	    \cline{1-2}
	    z & \multicolumn{2}{c}{} \\
	    \cline{1-1}
	\end{array}
	\right. \:
	\textrm{ for }x\leq y<z,
	\\
	\textrm{ and }
	\left.
	\begin{array}{c|c|c|}
	    \cline{3-3}
	     \multicolumn{1}{c}{} & & x \\
	    \cline{2-3}
	    & y & z  \\
	    \cline{2-3}
	\end{array}
	\equiv
	\begin{array}{|c|c|}
	    \cline{1-1}
	       x& \multicolumn{1}{c}{} \\
	    \cline{1-2}
	     y & z  \\
	    \cline{1-2}
	\end{array}	
	\equiv
	\begin{array}{|c|c|c}
	    \cline{1-2}
	     x &  z \\
	    \cline{1-2}
	    y & \multicolumn{2}{c}{} \\
	    \cline{1-1}
	\end{array}
	\right. \:	
	\textrm{ for }x<y\leq z	.
    \end{eqnarray*}
    
    \item  If $1<x<n$ and $x\leq y\leq\bar{x}$, then 
    \begin{eqnarray*}
	\left.
	\begin{array}{c|c|c|}
	    \cline{3-3}
	     \multicolumn{1}{c}{} & & y \\
	    \cline{2-3}
	    & x-1 & \overline{x-1}\rule{0pt}{2.5ex}  \\
	    \cline{2-3}
	\end{array}
	\equiv
	\begin{array}{c|c|c|}
	    \cline{2-3}
	      & x-1& y \\
	    \cline{2-3}
	    \multicolumn{1}{c}{}&  & \overline{x-1}\rule{0pt}{2.5ex}  \\
	    \cline{3-3}
	\end{array}	
	\equiv	
	\begin{array}{|c|c|c}
	    \cline{1-2}
	     x &  y \\
	    \cline{1-2}
	    \bar{x} & \multicolumn{2}{c}{} \\
	    \cline{1-1}
	\end{array}
	\right. \:
	\\
	\textrm{ and }
	\left.
	\begin{array}{c|c|c|}
	    \cline{3-3}
	     \multicolumn{1}{c}{} & & x \\
	    \cline{2-3}
	    & y & \overline{x}  \\
	    \cline{2-3}
	\end{array}
	\equiv
	\begin{array}{|c|c|}
	    \cline{1-1}
	       x-1& \multicolumn{1}{c}{} \\
	    \cline{1-2}
	     y & \overline{x-1}\rule{0pt}{2.5ex}  \\
	    \cline{1-2}
	\end{array}	
	\equiv	
	\begin{array}{|c|c|c}
	    \cline{1-2}
	     x-1 &  \overline{x-1}\rule{0pt}{2.5ex} \\
	    \cline{1-2}
	    y & \multicolumn{2}{c}{} \\
	    \cline{1-1}
	\end{array}
	\right. \qquad\quad	
	.
    \end{eqnarray*}
    
    \item  If $x\leq n-1$, then
    \begin{eqnarray*}
	&&\phantom{\textrm{ and }}
	\left\{
	\begin{array}{c}
	    \left.
	    \begin{array}{c|c|c|}
		\cline{3-3}
		 \multicolumn{1}{c}{} & & \bar{n} \\
		\cline{2-3}
		& n & \overline{x}  \\
		\cline{2-3}
	    \end{array}
	    \equiv
	    \begin{array}{|c|c|}
		\cline{1-1}
		   \bar{n} & \multicolumn{1}{c}{} \\
		\cline{1-2}
		 n & \bar{x}  \\
		\cline{1-2}
	    \end{array}	
	    \equiv	    
	    \begin{array}{|c|c|c}
		\cline{1-2}
		 \bar{n} &  \bar{x} \\
		\cline{1-2}
		n & \multicolumn{2}{c}{} \\
		\cline{1-1}
	    \end{array}
	    \right. \:
	    \\ 
	    \\
	    \left.
	    \begin{array}{c|c|c|}
		\cline{3-3}
		 \multicolumn{1}{c}{} & & n \\
		\cline{2-3}
		& \bar{n} & \overline{x}  \\
		\cline{2-3}
	    \end{array}
	    \equiv
	    \begin{array}{|c|c|}
		\cline{1-1}
		   n& \multicolumn{1}{c}{} \\
		\cline{1-2}
		 \bar{n} & \bar{x}  \\
		\cline{1-2}
	    \end{array}	
	    \equiv	    
	    \begin{array}{|c|c|c}
		\cline{1-2}
		 n &  \bar{x} \\
		\cline{1-2}
		\bar{n} & \multicolumn{2}{c}{} \\
		\cline{1-1}
	    \end{array}
	    \right. \:
	    \\
	\end{array}
	\right.
	\\&&
	\textrm{ and }
	\left\{
	\begin{array}{c}
	    \left.
	    \begin{array}{c|c|c|}
		\cline{3-3}
		 \multicolumn{1}{c}{} & & \bar{n} \\
		\cline{2-3}
		& x & n  \\
		\cline{2-3}
	    \end{array}
	    \equiv
	    \begin{array}{c|c|c|}
		\cline{2-3}
		  & x& \bar{n} \\
		\cline{2-3}
		\multicolumn{1}{c}{}&  & n  \\
		\cline{3-3}
	    \end{array}	
	    \equiv	    
	    \begin{array}{|c|c|c}
		\cline{1-2}
		 x &  \bar{n} \\
		\cline{1-2}
		n & \multicolumn{2}{c}{} \\
		\cline{1-1}
	    \end{array}
	    \right. \:
	    \\ 
	    \\
	    \left.
	    \begin{array}{c|c|c|}
		\cline{3-3}
		 \multicolumn{1}{c}{} & & n \\
		\cline{2-3}
		& x & \overline{n}  \\
		\cline{2-3}
	    \end{array}
	    \equiv
	    \begin{array}{c|c|c|}
		\cline{2-3}
		  & x& n \\
		\cline{2-3}
		\multicolumn{1}{c}{}&  & \bar{n}  \\
		\cline{3-3}
	    \end{array}	
	    \equiv	    
	    \begin{array}{|c|c|c}
		\cline{1-2}
		 x &  n \\
		\cline{1-2}
		\bar{n} & \multicolumn{2}{c}{} \\
		\cline{1-1}
	    \end{array}
	    \right. \:
	    \\
	\end{array}
	\right. .
    \end{eqnarray*}
    
    \item 
    \begin{eqnarray*}
	&&\phantom{\textrm{ and }}
	\left\{
	\begin{array}{c}
	    \left.
	    \begin{array}{c|c|c|}
		\cline{3-3}
		 \multicolumn{1}{c}{} & & n \\
		\cline{2-3}
		& \bar{n} & \bar{n}  \\
		\cline{2-3}
	    \end{array}
	    \equiv
	    \begin{array}{|c|c|}
		\cline{1-1}
		  n-1 & \multicolumn{1}{c}{} \\
		\cline{1-2}
		 \bar{n} & \overline{n-1}\rule{0pt}{2.5ex}  \\
		\cline{1-2}
	    \end{array}	
	    \equiv	    
	    \begin{array}{|c|c|c}
		\cline{1-2}
		 n-1 &  \overline{n-1}\rule{0pt}{2.5ex} \\
		\cline{1-2}
		\bar{n} & \multicolumn{2}{c}{} \\
		\cline{1-1}
	    \end{array}
	    \right. \:
	    \\ 
	    \\
	    \left.
	    \begin{array}{c|c|c|}
		\cline{3-3}
		 \multicolumn{1}{c}{} & & \bar{n} \\
		\cline{2-3}
		& n & n  \\
		\cline{2-3}
	    \end{array}
	    \equiv
	    \begin{array}{|c|c|}
		\cline{1-1}
		   n-1& \multicolumn{1}{c}{} \\
		\cline{1-2}
		 n & \overline{n-1}\rule{0pt}{2.5ex}  \\
		\cline{1-2}
	    \end{array}	
	    \equiv	    
	    \begin{array}{|c|c|c}
		\cline{1-2}
		 n-1 &  \overline{n-1}\rule{0pt}{2.5ex} \\
		\cline{1-2}
		n & \multicolumn{2}{c}{} \\
		\cline{1-1}
	    \end{array}
	    \right. \:
	    \\	    
	\end{array}
	\right.
	\\
	&&
	\textrm{ and }
	\left\{
	\begin{array}{c}
	    \left.
	    \begin{array}{c|c|c|}
		\cline{3-3}
		 \multicolumn{1}{c}{} & & \bar{n} \\
		\cline{2-3}
		& n-1 & \overline{n-1} \rule{0pt}{2.5ex} \\
		\cline{2-3}
	    \end{array}
	    \equiv
	    \begin{array}{c|c|c|}
		\cline{2-3}
		  & n-1& \bar{n} \\
		\cline{2-3}
		\multicolumn{1}{c}{}&  & \overline{n-1}\rule{0pt}{2.5ex}  \\
		\cline{3-3}
	    \end{array}	
	    \equiv	    
	    \begin{array}{|c|c|c}
		\cline{1-2}
		 \bar{n} &  \bar{n} \\
		\cline{1-2}
		n & \multicolumn{2}{c}{} \\
		\cline{1-1}
	    \end{array}
	    \right. \:
	    \\ 
	    \\
	    \left.
	    \begin{array}{c|c|c|}
		\cline{3-3}
		 \multicolumn{1}{c}{} & & n \\
		\cline{2-3}
		& n-1 & \overline{n-1}\rule{0pt}{2.5ex}  \\
		\cline{2-3}
	    \end{array}
	    \equiv
	    \begin{array}{c|c|c|}
		\cline{2-3}
		  & n-1& n \\
		\cline{2-3}
		\multicolumn{1}{c}{}&  & \overline{n-1}\rule{0pt}{2.5ex}  \\
		\cline{3-3}
	    \end{array}	
	    \equiv	    
	    \begin{array}{|c|c|c}
		\cline{1-2}
		 n &  n \\
		\cline{1-2}
		\bar{n} & \multicolumn{2}{c}{} \\
		\cline{1-1}
	    \end{array}
	    \right. \:
	    \\		        
	\end{array}
	\right. .
    \end{eqnarray*}
\end{enumerate}
If a word is composed entirely of barred letters or entirely of
unbarred letters, only relation ($1$) (the Knuth relation) applies,
and the type $A$ \textit{jeu de taquin} may be used.

\subsection{Properties of $B^{2,s}$}
As mentioned in the introduction, it was conjectured in~\cite{HKOTT:2001,HKOTY:1999}
that there are crystal bases $B^{r,s}$ associated with Kirillov--Reshetikhin 
modules $W^{r,s}$. In addition to the existence, Hatayama et al.~\cite{HKOTT:2001}
conjectured certain properties of $B^{r,s}$ which we state here in the specific
case of $B^{2,s}$ of type $D_n^{(1)}$.
\begin{Conj}[\cite{HKOTT:2001}]\label{conj:Brs}
If the crystal $B^{2,s}$ of type $D_n^{(1)}$ exists, it has the following
properties:
\begin{enumerate}
\item As a classical crystal $B^{2,s}$ decomposes as
$\bks{2}\cong \bigoplus_{k=0}^{s}B(k\Lt)$.
\item $B^{2,s}$ is perfect of level $s$.
\item $B^{2,s}$ is equipped with an energy function
$D_{B^{2,s}}$ such that $D_{B^{2,s}}(b)=k-s$ if $b$ is in the component of
$B(k\om_2)$ (in accordance with \eqref{eq:D single}).
\end{enumerate}
\end{Conj}

\section{Classical decomposition of $\tbs$}\label{sec:decomp}

In this section we begin our construction of the crystal $\tbs$
mentioned in Theorem \ref{thm:main}.  We do this by defining a
$U_q(D_n)$-crystal with vertices labelled by the set $\mathcal{T}(s)$
of tableaux of shape $(s,s)$ which satisfy conditions \ref{legal:1},
\ref{legal:2}, and \ref{legal:4} of Criterion \ref{legal}.
We will construct a bijection between $\mathcal{T}(s)$ and the
vertices of $\bigoplus_{i=0}^{s}B(i\Lt)$, so that $\mathcal{T}(s)$ may
be viewed as a $U_q(D_n)$-crystal with the classical decomposition
\eqref{eq:decomp}.  In section \ref{sec:aff} we will define $\ftil{0}$
and $\etil{0}$ on $\mathcal{T}(s)$ to give it the structure of a perfect
$U_q'(D_n^{(1)})$-crystal.  This crystal will be $\tbs$.

The reader may note in later sections that the main result of the paper 
does not depend on this explicit labeling of the vertices of $\tbs$.  
We have included it here because a description of the crystal in terms 
of tableaux will be needed to obtain a bijection with rigged configurations. 
It is through such a bijection that we anticipate being able to prove the $X=M$ 
conjecture for type $D$, as has already been done for special
cases in~\cite{OSS:2002,S:2004,SS:2004}.

\begin{Prop}\label{ablock}
Let $T\in\mathcal{T}(s)\setminus B(s\Lt)$ with $T\neq
\stack{1}{\bar{1}}\cdots \stack{1}{\bar{1}}$, and define
$\bar{\bar{i}}=i$ for $1\leq i \leq n$.  Then there is a unique
$a\in\{1,\ldots,n,\bar{n}\}$ and $m\in\mathbb{Z}_{>0}$ such that $T$
contains one of the following configurations (called an
$a$-configuration):
\begin{eqnarray*}
\stack{a}{b_{1}}
\underbrace{\stack{a}{\bar{a}} \cdots \stack{a}{\bar{a}}}_{m}
\stack{c_{1}}{d_{1}} ,&&
\textrm{where $b_{1}\neq \bar{a}$, and  
$c_{1}\neq a$ or $d_{1}\neq\bar{a}$};\\
\stack{b_{2}}{c_{2}}
\underbrace{\stack{a}{\bar{a}} \cdots \stack{a}{\bar{a}}}_{m}
\stack{d_{2}}{\bar{a}} ,&&
\textrm{where $d_{2}\neq a$, and  
$b_{2}\neq a$
or $c_{2}\neq\bar{a}$};\\
\stack{b_{3}}{c_{3}}
\underbrace{\stack{a}{\bar{a}} \cdots \stack{a}{\bar{a}}}_{m+1}
\stack{d_{3}}{e_{3}},&&
\textrm{where $b_{3}\neq a$ and $e_{3}\neq\bar{a}$}.
\end{eqnarray*}
\end{Prop}

\begin{proof}
If $s=1$, the set $\mathcal{T}(s)\setminus B(s\Lt)$ contains only
$\stack{1}{\bar{1}}$, so that the statement of the proposition is
empty.  Hence assume that $s\ge 2$.  The existence of an
$a$-configuration for some $a\in\{1,\ldots,n,\bar{n}\}$ follows from
the fact that $T$ violates condition~\ref{legal:3} of
Criterion~\ref{legal}.  The conditions on $b_{i}, c_{i}, d_{i}$ for
$i=1,2,3$ and $e_{3}$ can be viewed as stating that $m$ is chosen to
maximize the size of the $a$-configuration.  Condition~\ref{legal:1}
of Criterion~\ref{legal} and the conditions on the parameters
$b_i,c_i,d_i,e_3$ imply that there can be no other $a$-configurations
in $T$.
\end{proof}

The map $D_{2,s}:\mathcal{T}(s)\to \bigoplus_{k=0}^{s}B(k\Lt)$, called
the height-two drop map, is defined as follows for $T\in
\mathcal{T}(s)$.  If $T=\stack{1}{\bar{1}}\cdots \stack{1}{\bar{1}}$,
then $D_{2,s}(T)=\es\in B(0)$.  If $T\in B(s\Lt)$, $D_{2,s}(T)=T$.
Otherwise by Proposition~\ref{ablock}, $T$ contains a unique
$a$-configuration, and $D_{2,s}(T)$ is obtained from $T$ by removing
$\underbrace{\stack{a}{\bar{a}} \cdots \stack{a}{\bar{a}}}_{m}$.  

\begin{theorem}
    Let $T\in \mathcal{T}(s)$.  Then $D_{2,s}(T)$ satisfies Criterion
    \ref{legal}, and is therefore a tableau in
    $\bigoplus_{k=0}^{s}B(k\Lt)$.
\end{theorem}

\begin{proof}
    Condition $1$ is satisfied since the relation $\leq$ on our
    alphabet is transitive.  Conditions $2$ and $5$ are automatically
    satisfied, since the columns that remain are not changed.
    Condition $3$ is satisfied since by Proposition \ref{ablock},
    there can be no more than one $a$-configuration in $T$.  Condition
    $4$ is satisfied since $D_{2,s}$ does not remove any columns of
    the form $\stack{n-1}{n}$, $\stack{n-1}{\bar{n}}$,
    $\stack{n}{\overline{n-1}}$, or $\stack{\bar{n}}{\overline{n-1}}$.
\end{proof}

In Proposition \ref{filling}, we will show that $D_{2,s}$ is a
bijection by constructing its inverse.

\begin{Example}
We have
\begin{displaymath}
T=    
\left.
\begin{array}{|c|c|c|c|}
    \hline
    1 & 2 & 3& 3 \\
    \hline
    \overline{4} & \overline{2} & \overline{2} & \overline{1} \\
    \hline
\end{array}
\right. ,\quad
D_{2,4}(T)=
\left.
\begin{array}{|c|c|c|}
    \hline
    1 & 3& 3 \\
    \hline
    \overline{4} & \overline{2} & \overline{1} \\
    \hline
\end{array}
\right.  .
\end{displaymath}
\end{Example}

The inverse of $D_{2,s}$ is the height-two fill map
$F_{2,s}:\bigoplus_{k=0}^{s}B(k\Lt) \to\mathcal{T}(s)$.  Let
$t=\stack{a_{1}}{b_{1}}\cdots\stack{a_{k}}{b_{k}} \in B(k\Lt)$.  If
$k=s$, $F_{2,s}(t)=t$.  If $k<s$, then $F_{2,s}(t)$ is obtained by
finding a subtableau $\stack{a_{i}}{b_{i}}\stack{a_{i+1}}{b_{i+1}}$ in
$t$ such that
\begin{Crit}\label{eq:crit}
    $$b_{i}\leq\bar{a}_{i}\leq b_{i+1}\textrm{  or  }a_{i}\leq
    \bar{b}_{i+1}\leq a_{i+1}.$$
\end{Crit}
(Recall that $\bar{\bar{i}}=i$ for $i\in\{1,\ldots,n\}$.)  Note that
the first pair of inequalities imply that $a_{i}$ is unbarred, and the
second pair of inequalities imply that $b_{i+1}$ is barred.  We may
therefore insert between columns $i$ and $i+1$ of $t$ either the
configuration $\underbrace{\stack{a_{i}}{\bar{a}_{i}} \cdots
\stack{a_{i}}{\bar{a}_{i}}}_{s-k}$ or
$\underbrace{\stack{\bar{b}_{i+1}}{b_{i+1}} \cdots
\stack{\bar{b}_{i+1}}{b_{i+1}}}_{s-k}$, depending on which part of
Criterion \ref{eq:crit} is satisfied.  We say that $i$ is the filling
location of $t$.  If no such subtableau exists, then $F_{2,s}$ will
either append $\underbrace{\stack{a_{k}}{\bar{a}_{k}} \cdots
\stack{a_{k}}{\bar{a}_{k}}}_{s-k}$ to the end of $t$, or prepend
$\underbrace{\stack{\bar{b}_{1}}{b_{1}} \cdots
\stack{\bar{b}_{1}}{b_{1}}}_{s-k}$ to $t$.  In these cases the filling
locations are $k$ and $0$, respectively.

\begin{Prop} \label{filling}
The map $F_{2,s}$ is well-defined on $\bigoplus_{i=0}^{s}B(i\Lt)$.
\end{Prop}

The proof of this proposition follows from the next three lemmas.

\begin{Lemma} 
Suppose that $t\in\bigoplus_{k=0}^{s-1}B(k\Lt)$ has no subtableaux
$\stack{a_{i}}{b_{i}}\stack{a_{i+1}}{b_{i+1}}$ satisfying Criterion
\ref{eq:crit}.  Then exactly one of either appending
$\stack{a_{k}}{\bar{a}_{k}} \cdots \stack{a_{k}}{\bar{a}_{k}}$ or
prepending $\stack{\bar{b}_{1}}{b_{1}} \cdots
\stack{\bar{b}_{1}}{b_{1}}$ to $t$ will produce a tableau in
$\mathcal{T}(s)\setminus B(s\Lt)$.
\end{Lemma}

\begin{proof} Suppose
$t=\stack{a_{1}}{b_{1}}\cdots\stack{a_{k}}{b_{k}}\in B(k\Lt)$ is as
above for $k<s$.  We will show that if prepending
$\stack{\bar{b}_{1}}{b_{1}} \cdots \stack{\bar{b}_{1}}{b_{1}}$ to $t$
does not produce a tableau in $\mathcal{T}(s)\setminus B(s\Lt)$, then
appending $\stack{a_{k}}{\bar{a}_{k}} \cdots
\stack{a_{k}}{\bar{a}_{k}}$ to $t$ will produce a tableau in
$\mathcal{T}(s)\setminus B(s\Lt)$.  There are two reasons we might not
be able to prepend $\stack{\bar{b}_{1}}{b_{1}} \cdots
\stack{\bar{b}_{1}}{b_{1}}$; $b_{1}$ may be unbarred, or we may have
$a_{1}<\bar{b}_{1}$.

First, suppose $b_{1}$ is unbarred.  If $b_{k}$ is also unbarred, then
$b_{k}$ is certainly less than $\bar{a}_{k}$, so we may append 
$\stack{a_{k}}{\bar{a}_{k}} \cdots \stack{a_{k}}{\bar{a}_{k}}$ to $t$.  
Hence, suppose that $b_{k}$ is barred. We will show that $a_{k}$ is unbarred
and $\bar{a}_{k}>b_{k}$.

We know that $t$ has a subtableau of the form $\stack{a_{i}}{b_{i}}
\stack{a_{i+1}}{b_{i+1}}$ such that $b_{i}$ is unbarred and
$b_{i+1}$ is barred.  It follows that $a_{i}$ is unbarred, and
therefore $\bar{a}_{i}>b_{i}$.  Since $\stack{a_{i}}{b_{i}}
\stack{a_{i+1}}{b_{i+1}}$ does not satisfy Criterion~\ref{eq:crit}, this means 
that $\bar{a}_{i}>b_{i+1}$, which is equivalent to $\bar{b}_{i+1}>a_{i}$.  Once again 
observing that $\stack{a_{i}}{b_{i}}\stack{a_{i+1}}{b_{i+1}}$ does not satisfy 
Criterion~\ref{eq:crit}, this implies that $\bar{b}_{i+1}>a_{i+1}$; i.e., $a_{i+1}$ 
is unbarred, and $\bar{a}_{i+1}>b_{i+1}$.

We proceed with an inductive argument on $i<j<k$.  Suppose that 
$\stack{a_j}{b_j}\stack{a_{j+1}}{b_{j+1}}$ is a subtableau of $t$ such that $b_j$
and $b_{j+1}$ are barred, $a_j$ is unbarred, and $\bar{a}_j>b_j$.  By reasoning identical 
to the above, we conclude that
\begin{eqnarray}
     & \bar{a}_j>b_{j+1}
\Rightarrow
\bar{b}_{j+1}>a_j
\Rightarrow
\bar{b}_{j+1}>a_{j+1}
\Rightarrow
\bar{a}_{j+1}>b_{j+1}, &  \label{eq:imp}
\end{eqnarray}
which once again means that $a_{j+1}$ is unbarred.

This inductively shows that $a_{k}$ is unbarred and
$\bar{a}_{k}>b_{k}$, so we may append $\stack{a_{k}}{\bar{a}_{k}}
\cdots \stack{a_{k}}{\bar{a}_{k}}$ to $t$ to get a tableau in
$\mathcal{T}(s)\setminus B(s\Lt)$.  By a symmetrical argument, we
conclude that if $a_{k}$ is barred, then we may prepend
$\stack{\bar{b}_{1}}{b_{1}} \cdots \stack{\bar{b}_{1}}{b_{1}}$ to $t$.

Now, suppose that $b_{1}$ is barred and $\bar{b}_{1}>a_{1}$.  This
means that $a_{1}$ is unbarred and $\bar{a}_{1}>b_{1}$, so the
induction carried out in equation \ref{eq:imp} applies.  It follows that 
$a_{k}$ is unbarred and $\bar{a}_{k}>b_{k}$, so once again we may append 
$\stack{a_{k}}{\bar{a}_{k}} \cdots \stack{a_{k}}{\bar{a}_{k}}$ to $t$.  
Also, by a symmetrical argument, when $a_{k}$ is unbarred and $b_{k}>\bar{a}_{k}$, 
we may prepend $\stack{\bar{b}_{1}}{b_{1}} \cdots \stack{\bar{b}_{1}}{b_{1}}$ to $t$. 
Thus, when no subtableau of $t$ satisfy Criterion~\ref{eq:crit},
either appending $\stack{a_{k}}{\bar{a}_{k}} \cdots \stack{a_{k}}{\bar{a}_{k}}$ or 
prepending $\stack{\bar{b}_{1}}{b_{1}} \cdots \stack{\bar{b}_{1}}{b_{1}}$ to $t$ will 
produce a tableau in $\mathcal{T}(s)\setminus B(s\Lt)$.
\end{proof}

\begin{Lemma} 
Any tableau
$t=\stack{a_{1}}{b_{1}}\cdots\stack{a_{k}}{b_{k}}\in\bigoplus_{k=0}^{s-1}B(k\Lt)$
has no more than two filling locations.  If it has two, they are
consecutive integers, and this choice has no effect on $F_{2,s}(t)$.
\end{Lemma}

\begin{proof} Let $0\le i_*\le k$ be minimal such that $i_*$ is a
filling location of $t$.  First assume that $0<i_*<k$.  This implies
the existence of a subtableau
$\stack{a_{i_*}}{b_{i_*}}\stack{a_{i_*+1}}{b_{i_*+1}}$ which satisfies
Criterion~\ref{eq:crit}.

Suppose that the first condition $b_{i_{*}}\leq\bar{a}_{i_{*}}\leq
b_{i_{*}+1}$ of Criterion~\ref{eq:crit} is satisfied, and consider
whether $i_{*}+1$ can be a filling location.  If
$b_{i_{*}+1}\leq\bar{a}_{i_{*}+1}\leq b_{i_{*}+2}$, we have
$$b_{i_{*}+1}\leq\bar{a}_{i_{*}+1}\leq\bar{a}_{i_{*}}\leq
b_{i_{*}+1},$$
which implies that $\bar{a}_{i_{*}}=\bar{a}_{i_{*}+1}=b_{i_{*}+1}$, so
that $t$ violates part \ref{legal:3} of Criterion~\ref{legal}.
Similarly, if $a_{i_{*}+1}\leq\bar{b}_{i_{*}+2}\leq a_{i_{*}+2}$, then
we have $$\bar{a}_{i_{*}+1}\leq\bar{a}_{i_{*}}\leq b_{i_{*}+1}\leq
b_{i_{*}+2}\leq\bar{a}_{i_{*}+1},$$
which also implies that
$\bar{a}_{i_{*}}=\bar{a}_{i_{*}+1}=b_{i_{*}+1}$, once again violating
part~\ref{legal:3} of Criterion \ref{legal}.  We conclude that if
$i_{*}$ is a filling location for which Criterion \ref{eq:crit} is
satisfied by $b_{i_{*}}\leq\bar{a}_{i_{*}}\leq b_{i_{*}+1}$, then
$i_{*}+1$ is not a filling location.  Furthermore, this argument shows
that $a_{i_{*}+1}>a_{i_{*}}$ or $b_{i_{*}+1}>\bar{a}_{i_{*}}$.  By the
partial ordering on our alphabet, it follows that $t$ has no other
filling locations.

Now, suppose for the filling location $i_{*}$, Criterion \ref{eq:crit}
is satisfied by $a_{i_{*}}\leq\bar{b}_{i_{*}+1}\leq a_{i_{*}+1}$.  The
condition $a_{i_{*}+1}\leq\bar{b}_{i_{*}+2}\leq a_{i_{*}+2}$ for
$i_{*}+1$ to be a filling location implies that
$$\bar{b}_{i_{*}+2}\leq\bar{b}_{i_{*}+1}\leq
a_{i_{*}+1}\leq\bar{b}_{i_{*}+2},$$
which as above leads to a violation of part~\ref{legal:3} of Criterion
\ref{legal}.  However, $i_{*}+1$ may be a filling location if
Criterion \ref{eq:crit} is satisfied by
$b_{i_{*}+1}\leq\bar{a}_{i_{*}+1}\leq b_{i_{*}+2}$.  Note that this
inequality implies that $a_{i_{*}+1}\leq\bar{b}_{i_{*}+1}$, which
tells us that $a_{i_{*}+1}=\bar{b}_{i_{*}+1}$.  Thus, choosing to
insert
$\stack{\bar{b}_{i_{*}+1}}{b_{i_{*}+1}}\cdots\stack{\bar{b}_{i_{*}+1}}
{b_{i_{*}+1}}$ between columns $i_{*}$ and $i_{*}+1$ or to insert
$\stack{a_{i_{*}+1}}{\bar{a}_{i_{*}+1}}\cdots\stack{a_{i_{*}+1}}{\bar{a}_{i_{*}+1}}$
between columns $i_{*}+1$ and $i_{*}+2$ does not change $F_{2,s}(t)$.
Since $i_{*}+1$ is a filling location with Criterion~\ref{eq:crit}
satisfied by $b_{i_{*}}\leq\bar{a}_{i_{*}}\leq b_{i_{*}+1}$, the
preceding paragraph implies that there are no other filling locations
in $t$.

Finally, suppose that $i_*=0$ is a filling location for $t$; i.e.,
$b_{1}$ is barred, $a_{1}$ is unbarred, and $\bar{b}_{1}\leq a_{1}$.
If $1$ is a filling location, Criterion \ref{eq:crit} is satisfied by
$b_{1}\leq\bar{a}_{1}\leq b_{2}$; otherwise, part~\ref{legal:3} of
Criterion~\ref{legal} is violated.  Put together, this means that
$\bar{a}_{1}=b_{1}$, so prepending $\stack{\bar{b}_{1}}{b_{1}}\cdots
\stack{\bar{b}_{1}}{b_{1}}$ to $t$ and inserting
$\stack{a_{1}}{\bar{a}_{1}} \cdots\stack{a_{1}}{\bar{a}_{1}}$ between
columns $1$ and $2$ results in the same tableau.  As in the above
cases, part~\ref{legal:3} of Criterion \ref{legal} and the partial
order on the alphabet prohibit any other filling locations.
\end{proof}

\begin{Example}
Let $s=4$.  Then
\begin{displaymath}
t=
\left.
\begin{array}{|c|c|c|}
    \hline
    1 & 2& 3 \\
    \hline
    \overline{4} & \overline{2} & \overline{1} \\
    \hline
\end{array}
\right.  ,\quad
F_{2,4}(t)=
\left.
\begin{array}{|c|c|c|c|}
    \hline
    1 & 2 & 2& 3 \\
    \hline
    \overline{4} & \overline{2} & \overline{2} & \overline{1} \\
    \hline
\end{array}
\right. .
\end{displaymath}
While we could choose either column two or column three as the filling 
location, either choice results in the same tableau.
\end{Example}

\begin{Lemma} 
If a filling location of
$t=\stack{a_{1}}{b_{1}}\cdots\stack{a_{k}}{b_{k}}\in
\bigoplus_{i=0}^{s-1}B(i\Lt)$ satisfies Criterion \ref{eq:crit} with
both inequalities, then $F_{2,s}(t)$ is independent of this choice.
\end{Lemma}

\begin{proof} Suppose that $i_{*}\neq 0,k$ is a filling location for
$t$ where both parts of Criterion \ref{eq:crit} are satisfied.  This
means that the subtableau
$\stack{a_{i_{*}}}{b_{i_{*}}}\stack{a_{i_{*}+1}}{b_{i_{*}+1}}$
satisfies both $\bar{a}_{i_{*}}\leq b_{i_{*}+1}$ and $a_{i_{*}}\leq
\bar{b}_{i_{*}+1}$.  The latter of these implies that
$b_{i_{*}+1}\leq\bar{a}_{i_{*}}$, so we have
$\bar{a}_{i_{*}}=b_{i_{*}+1}$ and $\bar{b}_{i_{*}+1}=a_{i_{*}}$.
Thus, filling with either $\stack{a_{i_{*}}}{\bar{a}_{i_{*}}} \cdots
\stack{a_{i_{*}}}{\bar{a}_{i_{*}} }$ or
$\stack{\bar{b}_{i_{*}+1}}{b_{i_{*}+1}} \cdots
\stack{\bar{b}_{i_{*}+1}}{b_{i_{*}+1}}$ between columns $i_{*}$ and
$i_{*}+1$ results in the same tableau $F_{2,s}(t)$.
\end{proof}

\begin{Example}
To illustrate, for
\begin{displaymath}
t=
\left.
\begin{array}{|c|c|c|}
    \hline
    2 & 3& 3 \\
    \hline
    \overline{4} & \overline{2} & \overline{1} \\
    \hline
\end{array}
\right.  \quad \text{we have} \quad
F_{2,s}(t)=
\left.
\begin{array}{|c|c|c|c|}
    \hline
    2 & 2 & 3& 3 \\
    \hline
    \overline{4} & \overline{2} & \overline{2} & \overline{1} \\
    \hline
\end{array}
\right. .
\end{displaymath}
\end{Example}

By identifying $\mathcal{T}(s)$ with $\bigoplus_{i=0}^{s}B(i\Lt)$ via
the maps $D_{2,s}$ and $F_{2,s}$, we have defined a $U_q(D_n)$-crystal
with the decomposition \eqref{eq:decomp}, with vertices labelled by
the $2\times s$ tableaux of $\mathcal{T}(s)$.  The action of the
Kashiwara operators $\etil{i}$, $\ftil{i}$ for $i\in\{1,\ldots,n\}$ on
this crystal is defined in terms of the above bijection, given
explicitly by
\begin{equation} \label{eq:drop fill}
\begin{split}
    \etil{i}(T) &= F_{2,s}(\etil{i}(D_{2,s}(T)))\\
    \ftil{i}(T) &= F_{2,s}(\ftil{i}(D_{2,s}(T))), 
\end{split}
\end{equation}
for $T\in \mathcal{T}(s)$, where the $\etil{i}$ and $\ftil{i}$ on the
right are the standard Kashiwara operators on $U_q(D_n)$-crystals
\cite{KN:1994}.  In section~\ref{sec:aff} we will discuss the action of
$\etil{0}$ and $\ftil{0}$ on $\mathcal{T}(s)$, which will make $\mathcal{T}(s)$
into an affine crystal called $\tbs$.

Using the filling and dropping map we obtain a natural inclusion of 
$\mathcal{T}(s')$ into $\mathcal{T}(s)$ for $s'<s$.

\begin{Def} \label{def:upsilon}
For $s'<s$, the map 
$\Upsilon_{s'}^{s}:\mathcal{T}(s')\hookrightarrow\mathcal{T}(s)$ is 
defined by $\Upsilon_{s'}^{s}=F_{2,s}\circ D_{2,s'}$.
\end{Def}

\section{The branching component graph}\label{sec:bc}
The Dynkin diagram of $D_n^{(1)}$ has an automorphism interchanging
the nodes $0$ and $1$, which induces a map $\sigma:\bks{2}\to\bks{2}$ on the
crystals such that $\etil{0}=\sigma\etil{1}\sigma$ and $\ftil{0}=\sigma\ftil{1}\sigma$.
With this in mind, suppose we have defined $\ftil{0}$ on
$\mathcal{T}(s)$ to produce $\tbs$, and consider the following operations on
$\tbs$: Let $K\subset I$, and denote by $B_{K}$ the graph which
results from removing all $k$-colored edges from $\tbs$ for $k\in K$.
Then as directed graphs, we expect $B_{\{0\}}$ to be isomorphic to
$B_{\{1\}}$; otherwise, $\tbs$ and $\bks{2}$ will not be isomorphic.
We can gain some information about $\sigma$ by considering
$B_{\{0,1\}}$. The combinatorial structure of $B_{\{0,1\}}$ is encoded
in the branching component graph to be defined in this section.

The definition of $\sigma$ relies on several sets of data, which will be
defined in sections~\ref{sec:bc} and~\ref{sec:aff}. For all $k\geq 0$ there is 
a filtration of $B(k\om_2)$ by subgraphs isomorphic to $B(\ell\om_2)$ for 
$\ell\leq k$; this relates any classical component of $\tbs$ to the other classical 
components. Once this filtration is understood, we will see that the following 
data uniquely determines a vertex $b$ of $\tbs$:
\begin{enumerate}
\item its classical component $k$ in the direct sum $\bigoplus_{k=0}^sB(k\om_2)$;
\item its position $\ell$ in the filtration $B(k\om_2)\supset \cdots \supset 
B(\ell\om_2)\supset \cdots \supset B(0)$;
\item the number of $1$-arrows in a path to $b$ from the highest weight vector 
of $B(k\om_2)$;
\item the $D_{n-1}$-highest weight $\la$ of its connected component in $B_{\{0,1\}}$;
\item its position $b = \tilde{f}v_\la = \ftil{i_1}^{m_1}\ftil{i_2}^{m_2}\cdots v_\la$ 
in the $D_{n-1}$-crystal $B(\la)$.
\end{enumerate}	
The involution $\sigma$ has a very simple description in terms of these data.  
In fact, $\sigma$ changes only items (1) and (3), leaving the other data fixed.

\subsection{Definitions and preliminary discussion.}
The connected components of $B_{\{0,1\}}$ are $U_q(D_{n-1})$-crystals, indexed 
by partitions as described in this section. The decomposition of $\tbs$ into 
$B_{\{0,1\}}$ produces a branching component graph
for $\tbs$, which we denote $\Sup(\tbs)$.  
The vertices of this graph correspond to the connected
$U_q(D_{n-1})$-crystals; a vertex $v_{\la}$ is labelled
(non-uniquely) by the partition $\la$ indicating the classical highest
weight of the corresponding $U_q(D_{n-1})$-crystal.  The edges of
$\Sup(\tbs)$ are defined by placing an edge from $v_{\la}$ to
$v_{\mu}$ if there is a tableau $b\in B(v_{\la})$ such that
$\ftil{1}(b)\in B(v_{\mu})$, where $B(v_{\la})$ denotes the set of
tableaux contained in the $U_q(D_{n-1})$-crystal indexed by $v_{\la}$.

Note that it suffices to describe the decomposition of the component
of $\tbs$ with $U_{q}(D_{n})$ highest weight $k\Lt$ into
$U_q(D_{n-1})$-crystals for any $k\geq0$, since
\begin{equation*}
\Sup(\bigoplus_{k=0}^{s}B(k\Lt))=\bigoplus_{k=0}^{s}\Sup(B(k\Lt)).
\end{equation*}
Denote the branching component subgraph with classical highest weight
$k\Lt$ by $\Sup(k\Lt)$.  Since $\Sup(k\Lt)$ is determined by the
action of $\etil{i}$ and $\ftil{i}$ on $B(k\Lt)$ for $i=1,\ldots,n$,
which is in turn defined by composing the classical Kashiwara
operators with $D_{2,s}$ and $F_{2,s}$ (see equation~\eqref{eq:drop
fill}), it in fact suffices to determine the structure of
$\Sup(s\Lt)\subset\Sup(\tbs)$.

The branching component graph $\Sup(s\Lt)$ is characterized by the
following proposition.  We denote by $v_{s}$ the ``highest weight''
branching component vertex (that is to say the vertex $v$ such that
the highest weight vector $u_{s}$ of $B(s\Lt)$ is in $B(v)$) of
$\Sup(s\Lt)$.

\begin{Prop} \label{bc} The graph distance from $v_{s}$ defines a rank
function on $\Sup(s\Lt)$.  This graph has $2s+1$ ranks, and is
symmetric as a non-directed graph over rank $s$.  For $j\leq s$, the
$j^{\mathrm{th}}$ rank contains one of each partition
$\la=(\la_{1},\la_{2})\subset(s,j)$ such that $|\la|=s-j+2m$ for some
$m\in\mathbb{Z}_{\geq0}$.  For all ranks $0\leq j\leq 2s-1$, a vertex
$v_{\la}$ with rank $j$ has an arrow to a vertex $v_{\mu}$ with rank
$j+1$ if and only if $\la$ and $\mu$ are joined by an edge in Young's
lattice.
\end{Prop}

We begin by examining the first few ranks of $\Sup(s\Lt)$ in detail, 
then show that this proposition is true in general in 
sections~\ref{sec:content} and~\ref{sec:edge}.

The highest weight branching component vertex $v_{s}$ is indexed by the 
one-part partition $(s)$.  To see that this
is true, simply observe that the highest weight tableau of $B(s\Lt)$
is $\underbrace{ \stack{1}{2}\cdots \stack{1}{2}}_{s},$ and acting by
$\ftil{2},\ldots,\ftil{n}$ in the most general possible way will
affect only the bottom row.  When we map these bottom row subtableaux
component\-wise by $a\mapsto a-1$ and $\bar{a}\mapsto \overline{a-1}$
to tableaux of shape $(s)$, and apply the same map to the colors of the
arrows, this is clearly isomorphic to the $U_q(D_{n-1})$-crystal with
highest weight $s\om_1$.

Now, consider what can result from acting on a tableau
$T=\stack{a_1}{b_1}\cdots\stack{a_s}{b_s}$ in $B(v_{s})$ by
$\ftil{1}$.  Since $a_{1}=\cdots=a_{s}=1$, this will turn $a_{s}$ into
a $2$.  There are two cases to consider: if $b_{s}=\bar{2}$, this
results in a tableau with a configuration $\stack{2}{\bar{2}}$ at the
right end (note that $\ftil{i}$, $\etil{i}$ for $i=2,\ldots,n$ do not
act on this subtableau); otherwise, it is a tableau with
$a_{1}=\cdots=a_{s-1}=1$ where some element of $U_q(D_{n-1})$ can act
on the rightmost column.  In either case, we can act with
$\etil{2},\ldots,\etil{n}$ to find a $U_q(D_{n-1})$ highest weight
vector $T'=\stack{a'_{1}}{b'_{1}}\cdots \stack{a'_{s}}{b'_{s}}$, where
we have $b'_{1}=\cdots=b'_{s-1}=2$; in the first case, we have
$b'_{s}=\bar{2}$, in the other, we have $b'_{s}=3$.  Remove those
parts of these tableaux on which $\etil{i}$ and $\ftil{i}$ for
$i=2,\ldots,n$ do not act; in both cases, we remove
$a'_{1},\ldots,a'_{s-1}$, and in the first case we also remove the
$\stack{2}{\bar{2}}$ at the end.  We then have a skew tableau, which
when rectified by Lecouvey $D$ equivalence (or, since there are no
barred letters remaining, \textit{jeu de taquin}), is either the
tableau $2\cdots2$ of shape $(s-1)$, or the tableau of shape $(s,1)$
with $2$'s in the first row and a $3$ in the second.  We conclude that
there are two vertices of rank $1$ in $\Sup(s\Lt)$, corresponding to
the partitions $(s-1)$ and $(s,1)$.

Before we generalize this construction, we have a few technical
remarks.

The number of $1$-arrows in a minimal path in the crystal graph
between the highest weight tableau and a tableau $T$ is the
``$\alpha_{1}$-height'' of $T$.  Thus, the function
\begin{equation*}
r_{s}(v)=d(v,v_{s})=\min_{P(v,v_{s})}\{\textrm{number of edges in
}P(v,v_{s})\}
\end{equation*}
where $P(v,v_{s})$ is the set of all paths from $v$ to $v_{s}$ in
$\Sup(s\Lt)$, is a rank function on $\Sup(s\Lt)$.

\begin{Def} \label{def:null}
A null-configuration of size $k$ is
    \begin{eqnarray*}
    &
    \underbrace{
    \stack{1}{\bar{2}}
    \cdots
    \stack{1}{\bar{2}}
    \;
    \stack{2}{\bar{1}}
    \cdots
    \stack{2}{\bar{1}}
    }_{k}
    &
    \qquad\textrm{ if $k$ is even},
    \\
    &
    \underbrace{
    \stack{1}{\bar{2}}
    \cdots
    \stack{1}{\bar{2}}
    \;
    \stack{2}{\bar{2}}
    \;
    \stack{2}{\bar{1}}
    \cdots
    \stack{2}{\bar{1}}
    }_{k}
    &
    \qquad\textrm{ if $k$ is odd},
    \end{eqnarray*}
where the number of $1$'s equals the number of $\bar{1}$'s and the
number of $2$'s equals the number of $\bar{2}$'s.
\end{Def}

Null-configurations are named thus because $\etil{i}$ and $\ftil{i}$
for $i=2,\ldots,n$ send $T$ to $0$, where $T$ is the $2\times s$
tableau which is a null-configuration of size $s$.  Therefore, $T$ is
the basis vector for the trivial representation of $U_q(D_{n-1})$ in
$\Sup(s\Lt)$.  Put another way, inserting a null-configuration into a
tableau $T$ has no effect on $\vare_{i}(T)$ or $\varphi_{i}(T)$ for
$i=2,\ldots,n$.  This generalizes the phenomenon we observed in the
case of $\stack{2}{\bar{2}}$.

\subsection{Content of rank $j$.} \label{sec:content}
We now characterize the partitions occuring in any rank $0\le j\le s$
of the branching component graph. (Ranks greater than $s$ will be defined 
by the $\dual$-duality of the crystal as defined in section \ref{sec:dual crystal}.)  
We defer the discussion of the edges of the branching component graph to
section~\ref{sec:edge}.

Let $T\in B(s\Lt)$.  We wish to determine the vertex $v_{\la}$ of
$\Sup(s\Lt)$ for which $T\in B(v_{\la})$, and also to determine
$r_{s}(v_{\la})$.  As demonstrated for ranks $0$ and $1$ above,
determine the parts of $T$ on which $\etil{i}$ and $\ftil{i}$ for
$i=2,\ldots,n$ do not act: this will be a null-configuration of size
$r_{2}$ (possibly of size $0$), $r_{1}$ many $1$'s in the first row
before the null-configuration, and $r_{3}$ many $\bar{1}$'s in the
second row after the null-configuration.  We can extract from these
data the pair
\begin{equation} \label{eq:t}
 (t_{1},t_{2})=(r_{1}+r_{2},r_{2}+r_{3}), 
\end{equation}
where $t_{1},t_{2}\leq s$.
By observing the number of times $1$ appears in a sequence
$i_{1},\ldots,i_{p}$ such that the highest weight vector of $B(s\Lt)$
is $u_{s}=\etil{i_{1}}\ldots\etil{i_{p}}T$, it is easily seen that
$r_{s}(v_{\la})=s-t_{1}+t_{2}$.

Consider the set $\mathcal{J}$ of tableaux such that
$s-t_{1}+t_{2}=j\leq s$.  We wish to determine the partitions $\la$
such that $T\in\mathcal{J}$ are in a $U_q(D_{n-1})$-crystal with
highest weight specified by $\la$.  First, note that
$|\la|=2s-t_{1}-t_{2}$, since this is precisely the number of boxes
where $\etil{i}$ and $\ftil{i}$ for $i=2,\ldots,n$ act non-trivially.
It follows that $|\la|=s+j-2t_{2}$, so $|\la|\equiv s+j$ (mod 2), and
since $t_{2}$ ranges from $0$ to $j$, we have $s-j\leq|\la|\leq s+j$.
Based on the definition of $\etil{i}$ and $\ftil{i}$ given in section
\ref{sec:class crys}, it is clear that
other than the $t_{1}+t_{2}$ boxes with $1$'s, $\bar{1}$'s, and the
null-configuration, a $U_q(D_{n-1})$-highest weight tableau must have
only $2$'s and $3$'s.  We may remove the irrelevant $t_{1}+t_{2}$
boxes from $T$ resulting in a skew tableau $\Tred$.  All the letters
in $\Tred$ are unbarred, so the Lecouvey relations applied to
$w_{\Tred}$ yield the column word of the rectification of $\Tred$ (we
call this rectified tableau the completely reduced form of $T$), whose
shape has no more than two parts.  Let $\mathcal{I}\subset\mathcal{J}$
be the set of $U_q(D_{n-1})$-highest weight tableaux with specified
values for $t_{1}$ and $t_{2}$.  Then $\mathcal{I}$ includes tableaux
where the number of $2$'s ranges from $s-t_{2}$ up to
$\min(2s-t_{1}-t_{2}, s)$, and the number of $3$'s ranges
simultaneously from $s-t_{1}$ down to $\max(0, s-t_{1}-t_{2})$.  The
algorithm described above can therefore produce a tableau of any shape
$\la$ with two parts such that $|\la|=2s-t_{1}-t_{2}$, $\la_{1}\leq
s$, and $\la_{2}\leq s-t_{1}=j-t_{2}$.  By properties of the plactic
monoid, no two $U_q(D_{n-1})$-highest weight tableaux in $\mathcal{J}$
correspond to the same partition.

To summarize: In rank $j\leq s$ of $\Sup(s\Lt)$, the vertices
correspond exactly to partitions $\la=(\la_{1},\la_{2})\subset(s,j)$
such that $|\la|=s-j+2m$ for some $m\in\mathbb{Z}_{\geq0}$.

By the $\dual$-symmetry of $B(s\Lt)$ as described in 
section~\ref{sec:dual crystal}, it is clear that the
$U_q(D_{n-1})$-crystals of rank $j$ are the same as the
$U_q(D_{n-1})$-crystals of rank $2s-j$.  This completely characterizes
the vertices of $\Sup(s\Lt)$ by rank, and leads us to the following
remark.

\begin{Rem} If we consider the embedding
$U_q(D_{n-1})\hookrightarrow U_q(D_n)$ as implicitly described above,
and think of the action of $e_{1}, f_{1}\in U_q(D_n)$ as specifying a
rank function on the embedded $U_q(D_{n-1})$-modules in a given
$U_q(D_n)$-module with highest weight $s\Lt$, this provides a
combinatorial proof that the ranks are multiplicity-free.
\end{Rem}

\subsection{Edges of $\Sup(s\Lt)$.} \label{sec:edge}
We must now confirm that the pairs of vertices which have an arrow
between them are precisely those $v_{\la}$ and $v_{\mu}$ such that
$r_{s}(v_{\la})=j$ and $r_{s}(v_{\mu})=j+1$ for some $0\leq j\leq
2s-1$, and for which $\la$ and $\mu$ are adjacent in Young's lattice,
that is, $\mu$ is obtained from $\la$ by either adding or removing a
box.  To do this, we will construct tableaux in $B(v_{\la})$ such that
the shape of the completely reduced form of their image under
$\ftil{1}$ is the result of adding a box to $\la$.  The question of
removing boxes from $\la$ then is simply a matter of appealing to the
$\dual$-symmetry of the crystal graph as described in
section~\ref{sec:dual crystal}.

Our analysis breaks into two cases, where our tableau $T\in
B(v_{\la})$ may be of one of the following two forms:
\begin{enumerate}
    \item 
    $$
    T=
    \underbrace{
    \stack{1}{b_{1}}
    \cdots
    \cdots
    \stack{1}{b_{r_{1}}}}_{r_{1}}
    \underbrace{
    \stack{1}{\bar{2}}
    \cdots
    \stack{1}{\bar{2}}
    \bigg(
    \stack{2}{\bar{2}}
    \bigg)
    \stack{2}{\bar{1}}
    \cdots
    \stack{2}{\bar{1}}}_{r_{2}}
    \underbrace{
    \stack{a_{s-r_{3}+1}}{\bar{1}}
    \cdots
    \stack{a_{s}}{\bar{1}}}_{r_{3}} ,
    $$

    \item 
    $$
    T=
    \underbrace{
    \stack{1}{b_{1}}
    \cdots
    \stack{1}{b_{r_{1}}}}_{r_{1}}
    \underbrace{
    \stack{a_{r_{1}+1}}{b_{r_{1}+1}}
    \cdots
    \stack{a_{s-r_{3}}}{b_{s-r_{3}}}}_{u=s-r_{1}-r_{3}}
    \underbrace{
    \stack{a_{s-r_{3}+1}}{\bar{1}}
    \cdots
    \stack{a_{s}}{\bar{1}}}_{r_{3}} ,
    $$
\end{enumerate}
where in case ($1$), the block of length $r_{2}$ is a maximal
null-configuration, and in case ($2$), $a_{r_{1}+1}\neq 1$ and
$b_{s-r_{3}}\neq \bar{1}$ (we set $r_{2}=0$ here).  We now determine
for which partitions $\mu$ we can have $\ftil{1}(T)\in B(v_{\mu})$.
Recall from the previous subsection that for $T\in B(v_{\la})$, we
defined $\Tred$ to be the skew tableau which results from removing all
$1$'s, $\bar{1}$'s, and null-configurations from $T$.  Observe that
\begin{equation*}
    w_{\Tred}= \begin{cases}
    b_{1}\cdots b_{r_{1}}a_{s-r_{3}+1}\cdots a_{s} & \textrm{for 
    case (1)}\\
    b_{1}\cdots b_{r_{1}}b_{r_{1}+1}a_{r_{1}+1}\cdots 
    b_{s-r_{3}}a_{s-r_{3}}a_{s-r_{3}+1}\cdots a_{s} & \textrm{for case 
    (2)}.
    \end{cases}
\end{equation*}

In either case, if $b_{r_{1}}=\bar{2}$, the size of the
null-configuration in $\ftil{1}(T)$ is $r_{2}+1$, since in case ($1$)
$\ftil{1}$ acts on the middle of the null-configuration, and in case
($2$) $\ftil{1}$ acts on $a_{r_{1}}=1$.  It follows that
$w_{\ftil{1}(T)^{\#}}$ is simply $w_{\Tred}$ with the $\bar{2}$
contributed by $b_{r_{1}}$ removed.  If $b_{r_{1}}\neq\bar{2}$, we see
that $w_{\ftil{1}(T)^{\#}}$ is $w_{\Tred}$ with a $2$ inserted from
$a_{s-r_{3}}$ in case ($1$), and from $a_{r_{1}}$ in case ($2$).
Since we are currently concerned with adding boxes to $\la$, let us
assume that $b_{r_{1}}\neq\bar{2}$, and analyze how inserting a $2$ as
above affects the shape of the rectifications of these words.

Our augmented words are
\begin{equation}\label{eq:redword}
    w_{\ftil{1}(T)^{\#}}= \begin{cases}
     b_{1}\cdots b_{r_{1}}2\,a_{s-r_{3}+1}\cdots a_{s} 
     & \text{for case (1)}\\
     b_{1}\cdots b_{r_{1}}2\,b_{r_{1}+1}a_{r_{1}+1}\cdots 
     b_{s-r_{3}}a_{s-r_{3}}a_{s-r_{3}+1}\cdots a_{s}
     & \text{for case (2).}
\end{cases}
\end{equation} 
Recall that we have assumed that $b_{r_{1}}\neq \bar{2}$, which in
turn implies that all letters $b_{1},\ldots,b_{r_{1}}$ are strictly
less than $\bar{2}$.  Using relation ($1$) of Lecouvey type $D$
equivalence, we may therefore move the $2$ from position $a_{r_{1}}$
to the second position in the word.  This new word begins
$b_{1}2b_{2}$, with $b_{2}>2$.  Since we may view all the plactic
operations on this word as sliding moves, the subword $b_{2}\cdots
a_{s}$ can be rectified to give a tableau with no more than two rows.
Thus, all we have done is added one box to our shape.

We now show that this process can add a box to the top row of $\la$
unless $\la_{1}=s$, and it can add a box to the bottom row unless
$\la_{2}=\la_{1}$.  In the $U_q(D_{n-1})$-crystal $B(v_{\la})$, we
know that there is a $U_{q}(D_{n-1})$ highest weight tableau $T_{\la}$
of the form
\begin{equation*}
T_{\la}= \begin{cases}
\:\;
\underbrace{
\stack{1}{2}
\cdots
\stack{1}{2}}_{r_{1}-\la_{2}}
\underbrace{
\stack{1}{3}
\cdots
\stack{1}{3}}_{\la_{2}}
\underbrace{
\stack{1}{\bar{2}}
\cdots
\stack{1}{\bar{2}}
\bigg(
\stack{2}{\bar{2}}
\bigg)
\stack{2}{\bar{1}}
\cdots
\stack{2}{\bar{1}}}_{r_{2}}
\underbrace{
\stack{2}{\bar{1}}
\cdots
\stack{2}{\bar{1}}}_{r_{3}} 
& \text{for case (1)}\\
\underbrace{
\stack{1}{2}
\cdots
\stack{1}{2}}_{r_{1}-(\la_{2}-u)}
\underbrace{
\stack{1}{3}
\cdots
\stack{1}{3}}_{\la_{2}-u}
\underbrace{
\stack{2}{3}
\cdots
\stack{2}{3}}_{u}
\underbrace{
\stack{2}{\bar{1}}
\cdots
\stack{2}{\bar{1}}}_{r_{3}}
& \text{for case (2).}
\end{cases}
\end{equation*}
Note that in case ($1$) we have $\la_{2}\leq r_{3}$ and in case ($2$) 
we have $\la_{2}-u\leq r_{3}$; otherwise, acting by $\etil{2}$ can turn 
another $3$ into a $2$.

These tableaux yield the words
\begin{equation*}
w_{T_{\la}^{\#}}=\begin{cases}
\;\:\:\underbrace{2\cdots2}_{r_{1}-\la_{2}}
\underbrace{3\cdots3}_{\la_{2}}
\underbrace{2\cdots2}_{r_{3}} 
& \text{for case (1)}\\
\underbrace{2\cdots2}_{r_{1}-(\la_{2}-u)}
\underbrace{3\cdots3}_{\la_{2}-u}
\underbrace{32\cdots32}_{2u}
\underbrace{2\cdots2}_{r_{3}}
& \text{for case (2)}.
\end{cases}
\end{equation*}
The completely reduced form of these tableaux is a
two-row tableau with $r_{1}+r_{3}-\la_{2}$ $2$'s in the top row and
$\la_{2}$ $3$'s in the bottom row, or $2u+r_{1}+r_{3}-\la_{2}$ $2$'s
in the top row and $\la_{2}$ $3$'s in the bottom row, respectively.
It is easy to see that by adding a $2$ to $w_{\Tred}$ as in~\eqref{eq:redword}, 
we simply add a box containing a $2$ to the top row
of the completely reduced form of $T_{\la}$.  Note that this procedure
fails precisely when $T_{\la}$ can have no $2$'s added to it, in which
case there are $s$ $2$'s in $T_{\la}$, and thus $\la_{1}=s$.

Now suppose that $\la_{1}-\la_{2}>0$, so that adding a box to the
second row will produce a legal diagram.  Consider
$\tilde{T}_{\la}=\ftil{2}^{\la_{1}-\la_{2}}(T_{\la})$ (note that
$\la_{1}$ is the number of $2$'s in $T_{\la}$).  This tableau is in
$B(v_{\la})$, so its completely reduced form has shape $\la$, and we
see that
\begin{equation*}
w_{\tilde{T}_{\la}^{\#}}=\begin{cases}
\:\:\:
\underbrace{3\cdots3}_{|\la|-r_{3}}
\underbrace{2\cdots2}_{\la_{2}} 
\underbrace{3\cdots3}_{r_{3}-\la_{2}}
& \text{for case (1)}\\
\underbrace{3\cdots3}_{|\la|-r_{3}-2u}
\underbrace{32\cdots32}_{2u}
\underbrace{2\cdots2}_{\la_{2}-u}
\underbrace{3\cdots3}_{r_{3}-\la_{2}+u}
& \text{for case (2).}
\end{cases}
\end{equation*}
The rectified tableau has $\la_{2}$ $2$'s followed by
$\la_{1}-\la_{2}$ $3$'s in the top row, and $\la_{2}$ $3$'s in the
bottom row.  From this description, we see that adding a $2$ to
$w_{\tilde{T}_{\la}^{\#}}$ as in \eqref{eq:redword} affects the
completely reduced tableau by preventing one of the $3$'s from the
bottom row from being slid up to the top row; i.e., $\la_{2}$ is
increased by $1$.  Since we add only one box at a time and the only
shape in rank $0$ is $(s,0)$, we know that the number of boxes in the
second row can never exceed the rank.

We now invoke the $\dual$-duality of the crystal graph to deal with
how boxes can be removed from $\la$.  If $v_{\la}\in\Sup(s\Lt)$ has
rank $p$, there is a unique vertex $v'_{\la}$, called the
complementary vertex of $v_{\la}$, with rank $2s-p$ for which the
corresponding $U_q(D_{n-1})$-crystal is $B(\la)$.  This involution
agrees with the $\dual$-crystal involution of section~\ref{sec:dual
crystal}.  We wish to show that there is an arrow from $v_{\la}$ to
$v_{\mu}$, where $\la/\mu$ is a single box and
$r_{s}(v_{\mu})=r_{s}(v_{\la})+1$.  Recall that by definition, this is
the case when for some $T\in B(v_{\la})$ we have $\ftil{1}(T)\in
B(v_{\mu})$.  Observe that $r_{s}(v'_{\la})=r_{s}(v'_{\mu})+1$, and
$\la$ is the result of adding a box to $\mu$; therefore, there is an
arrow from $v'_{\mu}$ to $v'_{\la}$.  It follows that we can find some
$T\in B(v'_{\mu})$ such that $\ftil{1}(T)\in B(v'_{\la})$.  In turn,
we have $T^{\dual}\in B(v_{\mu})$ and
$(\ftil{1}(T))^{\dual}=\etil{1}(T^{\dual})\in B(v_{\la})$.  Since we
know that $\ftil{1}(\etil{1}(T^{\dual}))=T^{\dual}\in B(v_{\mu})$, we
have shown that there is an arrow from $v_{\la}$ to $v_{\mu}$.

The arguments of sections~\ref{sec:content} and~\ref{sec:edge} prove
Proposition \ref{bc}.

\subsection{Construction of $\Sup(\tbs)$.}
Observe that $\Sup(\tbs)=\bigcup_{i=0}^{s}\Sup(i\Lt)$.  Let
$v_{\la}\in\Sup(i\Lt)\subset\Sup(\tbs)$.  Define
$R(v_{\la})=r_{i}(v_{\la})+s-i$.  This defines a rank on all of
$\Sup(\tbs)$.  Note that
$\Sup(i\Lt)\subset\Sup((i+1)\Lt)$, and this inclusion
is compatible with $R$.  Also note that if $R(v_{\la})=p$, then
$v'_{\la}$, the complementary vertex to $v_{\la}$, is now defined to
be the vertex of rank $2s-p$ with the same shape and in the same
component as $v_{\la}$.

To illustrate, $\Sup(\tilde{B}^{2,2})$ is given in Figure \ref{fig:sup}, with 
rank $0$ in the first line, rank $1$ in the second, etc.

\begin{figure}
\begin{center}\setlength{\unitlength}{1in}
    \begin{picture}(-2.2,0)
	\put(-2.5,0){\line(1,0){.25}}
	\put(-2.5,0){\line(0,-1){.125}}
	\put(-2.25,0){\line(0,-1){.125}}
	\put(-2.5,-.125){\line(1,0){.25}}
	\put(-2.375,0){\line(0,-1){.125}}
	\put(-2.45,-.16){\vector(-1,-2){.1}}
	\put(-2.3,-.16){\vector(1,-2){.12}}
	\put(-2.75,-.38){\line(1,0){.25}}
	\put(-2.75,-.38){\line(0,-1){.25}}
	\put(-2.625,-.38){\line(0,-1){.25}}
	\put(-2.75,-.505){\line(1,0){.25}}
	\put(-2.5,-.38){\line(0,-1){.125}}
	\put(-2.75,-.63){\line(1,0){.125}}
	\put(-2.85,-.68){\vector(-1,-2){.1}}
	\put(-2.65,-.68){\vector(0,-1){.2}}
	\put(-2.52,-.63){\vector(1,-1){.25}}
	\put(-2.25,-.44){\line(1,0){.125}}
	\put(-2.25,-.44){\line(0,-1){.125}}
	\put(-2.125,-.44){\line(0,-1){.125}}
	\put(-2.25,-.565){\line(1,0){.125}}
 	\put(-2.09,-.68){\vector(1,-2){.1}}
 	\put(-2.2,-.68){\vector(0,-1){.2}}
 	\put(-2.35,-.63){\vector(-1,-1){.25}}
	\put(-3.1,-.9){\line(1,0){.25}}
	\put(-3.1,-.9){\line(0,-1){.25}}
	\put(-2.975,-.9){\line(0,-1){.25}}
	\put(-3.1,-1.025){\line(1,0){.25}}
	\put(-2.85,-.9){\line(0,-1){.25}}
	\put(-3.1,-1.15){\line(1,0){.25}}
	\put(-2.95,-1.2){\vector(1,-2){.1}}
	\put(-2.75,-.96){\line(1,0){.25}}
	\put(-2.75,-.96){\line(0,-1){.125}}
	\put(-2.5,-.96){\line(0,-1){.125}}
	\put(-2.75,-1.085){\line(1,0){.25}}
	\put(-2.625,-.96){\line(0,-1){.125}}
	\put(-2.65,-1.2){\vector(0,-1){.2}}
 	\put(-2.55,-1.2){\vector(1,-1){.25}}
	\put(-2.3,-.9){\line(1,0){.125}}
	\put(-2.3,-.9){\line(0,-1){.25}}
	\put(-2.175,-.9){\line(0,-1){.25}}
	\put(-2.3,-1.025){\line(1,0){.125}}
	\put(-2.3,-1.15){\line(1,0){.125}}
	\put(-2.27,-1.2){\vector(-1,-1){.21}}
 	\put(-2.2,-1.2){\vector(0,-1){.2}}
	\put(-2,-1.08){$\emptyset$}
 	\put(-2,-1.2){\vector(-1,-2){.1}}
	\put(-2.75,-1.42){\line(1,0){.25}}
	\put(-2.75,-1.42){\line(0,-1){.25}}
	\put(-2.625,-1.42){\line(0,-1){.25}}
	\put(-2.75,-1.545){\line(1,0){.25}}
	\put(-2.5,-1.42){\line(0,-1){.125}}
	\put(-2.75,-1.67){\line(1,0){.125}}
	\put(-2.55,-1.7){\vector(1,-2){.1}}
	\put(-2.25,-1.48){\line(1,0){.125}}
	\put(-2.25,-1.48){\line(0,-1){.125}}
	\put(-2.125,-1.48){\line(0,-1){.125}}
	\put(-2.25,-1.608){\line(1,0){.125}}
	\put(-2.2,-1.66){\vector(-1,-2){.12}}	
	\put(-2.5,-1.95){\line(1,0){.25}}
	\put(-2.5,-1.95){\line(0,-1){.125}}
	\put(-2.25,-1.95){\line(0,-1){.125}}
	\put(-2.5,-2.075){\line(1,0){.25}}
	\put(-2.375,-1.95){\line(0,-1){.125}}
    \end{picture}
\end{center}
\begin{center}\setlength{\unitlength}{1in}
    \begin{picture}(-5.1,.45)
	\put(-2.25,.06){\line(1,0){.125}}
	\put(-2.25,.06){\line(0,-1){.125}}
	\put(-2.125,.06){\line(0,-1){.125}}
	\put(-2.25,-.065){\line(1,0){.125}}
	\put(-2.125,-.12){\vector(1,-2){.12}}
	\put(-2.25,-.12){\vector(-1,-2){.1}}
	\put(-2.42,-.4){\line(1,0){.125}}
	\put(-2.42,-.4){\line(0,-1){.25}}
	\put(-2.295,-.4){\line(0,-1){.25}}
	\put(-2.42,-.525){\line(1,0){.125}}
	\put(-2.42,-.65){\line(1,0){.125}}
	\put(-2.35,-.73){\vector(1,-2){.1}}
	\put(-2,-.58){$\emptyset$}
	\put(-2,-.68){\vector(-1,-2){.127}}	
	\put(-2.25,-.98){\line(1,0){.125}}
	\put(-2.25,-.98){\line(0,-1){.125}}
	\put(-2.125,-.98){\line(0,-1){.125}}
	\put(-2.25,-1.105){\line(1,0){.125}}
	\put(-0.9,-.58){$\emptyset$} 
   \end{picture}
\end{center}
\vspace{4cm}
\caption{Branching component graph $\Sup(\tilde{B}^{2,2})$\label{fig:sup}}
\end{figure}

\section{Affine Kashiwara operators}\label{sec:aff}
Since we know that $B_{\{0\}}$ and $B_{\{1\}}$ are isomorphic as
directed graphs, it is clear that we can put $0$-colored edges in the
branching component graph in such a way that interchanging the
$1$-edges and the $0$-edges and applying some shape-preserving
bijection $\hsig$ to the vertices of the branching component graph will 
produce an isomorphic colored
directed graph.  Such a bijection can be naturally extended to
$\sigma:\tbs\to\tbs$ as follows.  Let $b\in
B(v_{\la})\subset\tbs$ for some branching component vertex
$v_{\la}$, and let $u_{\la}$ denote the $U_q(D_{n-1})$-highest weight
vector of $B(v_{\la})$.  Then for some finite sequence
$i_{1},\ldots,i_{k}$ of integers in $\{2,\ldots,n\}$, we know that
$\ftil{i_{1}}\cdots\ftil{i_{k}}u_{\la}=b$.  Let
$v^{\dag}_{\la}=\hsig(v_{\la})$, and let $u^{\dag}_{\la}$ be the
highest weight vector of $B(v^{\dag}_{\la})$.  We may define
$\sigma(b)=\ftil{i_{1}}\cdots\ftil{i_{k}}u^{\dag}_{\la}$.  This
involution of $\tbs$ allows us to define the affine structure of 
the crystal by the following equations:
\begin{equation}    \label{eq:zeros}
	\ftil{0}=\sigma\ftil{1}\sigma\qquad\textrm{ and }\qquad
	\etil{0}=\sigma\etil{1}\sigma.
\end{equation}

\begin{Def} \label{def:B tilde}
    The affine crystal $\tbs$ is given by the set $\mathcal{T}(s)$ as 
    defined in section \ref{sec:decomp} with $\etil{i}, \ftil{i}$ 
    for $1\leq i\leq n$ as in \eqref{eq:drop fill} and $\etil{0}, 
    \ftil{0}$ as in \eqref{eq:zeros}.
\end{Def}

Using $\hsig$ as defined in section \ref{sec:hsig}, it will be
shown in section~\ref{sec:perfect} that the resulting
$U_q'(D_{n}^{(1)})$-crystal $\tbs$ is perfect.

\subsection{Construction of $\hsig$.}\label{sec:hsig}
We will define $\hsig(v_{\la})$ for $R(v_{\la})\leq s$, and observe
that $\hsig(v'_{\la})=\hsig(v_{\la})'$, where $v'$ denotes the
complementary vertex of $v$.  Let $v_{\la}\in\Sup(k\Lt)$,
$R(v_{\la})=p$, and $\ell$ be minimal such that
$\iotil_{k}^{s}(v_{\la})\in \iotil_{\ell}^{s}(\Sup(\ell\Lt))$, where
$\iotil_{i}^{j}$ is the embedding of $\Sup(i\Lt)$ in
$\Sup(j\Lt)$ for $i<j$.  Then by the inclusion
$\Sup(i\Lt)\subset\Sup((i+1)\Lt)$ for
$i=0,\ldots,s-1$, there are $s-\ell+1$ vertices of the same shape as
$v_{\la}$ of rank $p$ in $\Sup(\tbs)$, one in each $\Sup(j\Lt)$ for
$j=\ell,\ldots,s$.  We define $\hsig(v_{\la})$ to be the vertex of the
same shape as $v_{\la}$ of rank $2s-p$ in the component
$\Sup((s+\ell-k)\Lt)$.

The action of $\hsig$ on $\Sup(\tbt)$ is given in 
Figure \ref{fig:sigma}.  

\begin{figure}
\begin{center}\setlength{\unitlength}{1in}
    \begin{picture}(-2.2,0)
	\put(-2.5,0){\line(1,0){.25}}
	\put(-2.5,0){\line(0,-1){.125}}
	\put(-2.25,0){\line(0,-1){.125}}
	\put(-2.5,-.125){\line(1,0){.25}}
	\put(-2.375,0){\line(0,-1){.125}}
%

	\put(-2.75,-.38){\line(1,0){.25}}
	\put(-2.75,-.38){\line(0,-1){.25}}
	\put(-2.625,-.38){\line(0,-1){.25}}
	\put(-2.75,-.505){\line(1,0){.25}}
	\put(-2.5,-.38){\line(0,-1){.125}}
	\put(-2.75,-.63){\line(1,0){.125}}
%

	\put(-2.25,-.44){\line(1,0){.125}}
	\put(-2.25,-.44){\line(0,-1){.125}}
	\put(-2.125,-.44){\line(0,-1){.125}}
	\put(-2.25,-.565){\line(1,0){.125}}
%

	\put(-3.1,-.9){\line(1,0){.25}}
	\put(-3.1,-.9){\line(0,-1){.25}}
	\put(-2.975,-.9){\line(0,-1){.25}}
	\put(-3.1,-1.025){\line(1,0){.25}}
	\put(-2.85,-.9){\line(0,-1){.25}}
	\put(-3.1,-1.15){\line(1,0){.25}}
%

	\put(-2.75,-.96){\line(1,0){.25}}
	\put(-2.75,-.96){\line(0,-1){.125}}
	\put(-2.5,-.96){\line(0,-1){.125}}
	\put(-2.75,-1.085){\line(1,0){.25}}
	\put(-2.625,-.96){\line(0,-1){.125}}
%

	\put(-2.3,-.9){\line(1,0){.125}}
	\put(-2.3,-.9){\line(0,-1){.25}}
	\put(-2.175,-.9){\line(0,-1){.25}}
	\put(-2.3,-1.025){\line(1,0){.125}}
	\put(-2.3,-1.15){\line(1,0){.125}}
%

	\put(-2,-1.08){$\emptyset$}

	\put(-2.75,-1.42){\line(1,0){.25}}
	\put(-2.75,-1.42){\line(0,-1){.25}}
	\put(-2.625,-1.42){\line(0,-1){.25}}
	\put(-2.75,-1.545){\line(1,0){.25}}
	\put(-2.5,-1.42){\line(0,-1){.125}}
	\put(-2.75,-1.67){\line(1,0){.125}}
%

	\put(-2.25,-1.48){\line(1,0){.125}}
	\put(-2.25,-1.48){\line(0,-1){.125}}
	\put(-2.125,-1.48){\line(0,-1){.125}}
	\put(-2.25,-1.608){\line(1,0){.125}}
%
	
	\put(-2.5,-1.95){\line(1,0){.25}}
	\put(-2.5,-1.95){\line(0,-1){.125}}
	\put(-2.25,-1.95){\line(0,-1){.125}}
	\put(-2.5,-2.075){\line(1,0){.25}}
	\put(-2.375,-1.95){\line(0,-1){.125}}
    \end{picture}
\end{center}
\begin{center}\setlength{\unitlength}{1in}
    \begin{picture}(-5.1,.45)
	\put(-2.25,.06){\line(1,0){.125}}
	\put(-2.25,.06){\line(0,-1){.125}}
	\put(-2.125,.06){\line(0,-1){.125}}
	\put(-2.25,-.065){\line(1,0){.125}}
%

	\put(-2.42,-.4){\line(1,0){.125}}
	\put(-2.42,-.4){\line(0,-1){.25}}
	\put(-2.295,-.4){\line(0,-1){.25}}
	\put(-2.42,-.525){\line(1,0){.125}}
	\put(-2.42,-.65){\line(1,0){.125}}
%

	\put(-2,-.58){$\emptyset$}
%

	\put(-2.25,-.98){\line(1,0){.125}}
	\put(-2.25,-.98){\line(0,-1){.125}}
	\put(-2.125,-.98){\line(0,-1){.125}}
	\put(-2.25,-1.105){\line(1,0){.125}}
	\put(-0.9,-.58){$\emptyset$} 
%
%
	\qbezier[2000](-0.95,-.45)(-2.125,-0.2)(-3.36,-.45)
	\put(-0.94,-.45){\vector(4,-1){0}}
	\put(-3.37,-.45){\vector(-4,-1){0}}
	\put(-2.3,-1){\vector(-4,3){1.2}}
	\put(-2.3,-1){\vector(4,-3){0}}
	\put(-2.3,-.1){\vector(-4,-3){1.2}}
	\put(-2.3,-.1){\vector(4,3){0}}
	\qbezier[200](-2.46,-.65)(-3,-0.85)(-3.58,-.65)
	\put(-2.44,-.64){\vector(4,1){0}}
	\put(-3.6,-.64){\vector(-4,1){0}}
	\put(-3.85,.3){\vector(0,-1){1.72}}
	\put(-3.85,.32){\vector(0,1){0}}
	\qbezier[150](-4.15,-.18)(-4.324,-.53)(-4.15,-.88)
	\put(-4.14,-.16){\vector(1,2){0}}
	\put(-4.14,-.9){\vector(1,-2){0}}

   \end{picture}
\end{center}
\vspace{4cm} 
\caption{Definition of $\hsig$ on $\Sup(\tbt)$\label{fig:sigma}}
\end{figure}

\subsection{Combinatorial construction of $\sigma$.}

We can also give a direct combinatorial description of $\sigma(T)$ for
any $T\in\tbs$.  As an auxilliary construction (which will also be
useful in its own right later on), we combinatorially describe
$\iota_{i}^{j}:B(i\Lt)\hookrightarrow B(j\Lt)$, the unique crystal
embedding that agrees with
$\iotil_{i}^{j}:\Sup(i\Lt)\hookrightarrow\Sup(j\Lt)$.

\begin{Rem}
It will often be useful to identify $B(k\Lt)$ with 
its image in $\tbs$.  We will use the notation $T\in 
B(k\Lt)\subset\tbs$ to indicate this identification.
\end{Rem}

Let $i\in\{0,\ldots,s-1\}$, so $\iota_{i}^{i+1}$ denotes the embedding
of $B(i\Lt)$ in $B((i+1)\Lt)$.  Let $T\in B(i\Lt)\subset \tbs$.  This embedding
can be combinatorially understood through the following observations:
\begin{Rem} \label{rem:comb iota}
    \begin{tabular}{c}
    \end{tabular}
    \begin{itemize}
	\item $\varphi_{k}(T)=\varphi_{k}(\iota_{i}^{i+1}(T))$ and
	$\vare_{k}(T)=\vare_{k}(\iota_{i}^{i+1}(T))$ for $k=2,\ldots,n$;
    
	\item $D_{2,s}(\iota_{i}^{i+1}(T))$ has one more column than 
	$D_{2,s}(T)$ (recall $D_{2,s}$ from section \ref{sec:decomp});
	
	\item Let $v(T)$ be the branching component vertex containing $T$.  
	Then $R(v(T))=R(v(\iota_{i}^{i+1}(T)))$, so the rank of 
	$v(\iota_{i}^{i+1}(T))$ in $\Sup((i+1)\Lt)$ is one greater than the 
	rank of $v(T)$ in $\Sup(i\Lt)$.
    \end{itemize}
\end{Rem}
In other words, we know that $T$ has a maximal $a$-configuration of
size $s-i$ (section \ref{sec:decomp}), and has completely reduced form
$\Tred$ (section \ref{sec:bc}).  Furthermore, let $r_{1T}$ be the number
of $1$'s in the first row of $D_{2,s}(T)$ to the left of a
null-configuration, and similarly define $r_{2T}, r_{3T}, t_{1T}$, and
$t_{2T}$ as in \eqref{eq:t}.  Then the rank of $v(T)$ in $\Sup(i\Lt)$
is $i-t_{1T}+t_{2T}$.  We wish to construct a tableau $S$ with an
$a$-configuration of size $s-i-1$ such that $S^{\#}=\Tred$ and
$(i+1)-t_{1S}+t_{2S}=(i-t_{1T}+t_{2T})+1$; i.e.,
$t_{1S}-t_{2S}=t_{1T}-t_{2T}$.  Based on properties of the height $2$
type $D$ sliding algorithm of section \ref{sec:plac}, these conditions
can only be satisfied when $t_{jS}=t_{jT}+1$ for $j=1,2$.

We can calculate $\iota_{i}^{i+1}(T)$ by the following algorithm:
\begin{algorithm}\label{alg:iota}\mbox{}
\begin{enumerate}
    \item Remove the $a$-configuration of size $s-i$ from $T$ and
    slide it to get a $2\times i$ tableau.

    \item Remove the $1$'s, $\bar{1}$'s and the null-configuration 
    from the result to get a skew tableau of shape 
    $(i,i-t_{2T})/(t_{1T})$.
    
    \item Using the type $D$ sliding algorithm, produce a skew tableau
    of shape $((i+1),(i+1)-(t_{2T}+1))/(t_{1T}+1)$.

    \item Fill this tableau with $1$'s, $\bar{1}$'s, and a 
    null-configuration so that the result is a $2\times (i+1)$ tableau.
    
    \item Use the height $2$ fill map $F_{2,s}$ (section \ref{sec:decomp}) to 
    insert $s-i-1$ columns into the tableau.
\end{enumerate}
\end{algorithm}
This produces the unique tableau satisfying the three properties of
Remark~\ref{rem:comb iota}.

\begin{Example}
Let 
$$
T=
\begin{array}{|c|c|c|c|c|c|c|}
    \hline
    1 & 1 & 2 & 2 & 2 & \bar{3} & \bar{2}  \\
    \hline
    2 & 2 & 3 & \bar{2} & \bar{2} & \bar{2} & \bar{1}  \\
    \hline
     
\end{array}
\in B(5\Lt)\subset \tilde{B}^{2,7}.
$$
Running through the steps of our algorithm (using relation (2) of
section \ref{sec:plac} for step (3)) gives us
\begin{enumerate}
    \item 
    $
    \phantom{\iota_{5}^{6}(T)=}\:        \phantom{\iota(T)=}   
    \begin{array}{|c|c|c|c|c|}
	\hline
	1 & 1 & 2 & \bar{3} & \bar{2}  \\
	\hline
	2 & 2 & 3 & \bar{2} & \bar{1}  \\
	\hline
     
    \end{array}
    $
    \vspace{.1in}

    \item 
    $
    \phantom{\iota_{5}^{6}(T)=}\:    \phantom{\iota(T)=}
    \begin{array}{|c|c|c|c|c|c}
	\cline{3-5}
	 \multicolumn{1}{c}{} & & 2 & \bar{3} & \bar{2}  \\
	\hline
	2 & 2 & 3 & \bar{2} &  \multicolumn{1}{c}{} \\
	\cline{1-4}
     
    \end{array}
    $
    \vspace{.1in}

    \item 
    $
    \phantom{\iota_{5}^{6}(T)=}\:	\phantom{\iota(T)=}
    \begin{array}{|c|c|c|c|c|c|}
	\cline{4-6}
	 \multicolumn{2}{c}{}  &  & 3 & \bar{3} & \bar{2}  \\
	\hline
	2 & 2 & 3 & \bar{3} &  \multicolumn{2}{c}{}  \\
	\cline{1-4}
     
    \end{array}
    $    
    \vspace{.1in}
    
    \item 
    $
    \phantom{\iota_{5}^{6}(T)=}\:	\phantom{\iota(T)=}
    \begin{array}{|c|c|c|c|c|c|}
	\hline
	1 & 1 & 1 & 3 & \bar{3} & \bar{2}  \\
	\hline
	2 & 2 & 3 & \bar{3} & \bar{1} & \bar{1} \\
	\hline
     
    \end{array}
    $    
    \vspace{.1in}
    
    \item 
    $
    \phantom{\iota(T)=}
    \iota_{5}^{6}(T)=
    \begin{array}{|c|c|c|c|c|c|c|}
	\hline
	1 & 1 & 1 & 3 & 3 & \bar{3} & \bar{2}  \\
	\hline
	2 & 2 & 3 & \bar{3} & \bar{3} & \bar{1} & \bar{1} \\
	\hline
     
    \end{array}
    \in B(6\Lt)\subset \tilde{B}^{2,7}.
    $ 
\end{enumerate}
\end{Example}

Composing these maps gives us the following algorithm for 
calculating $\iota_{i}^{j}(T)$, where $T\in B(i\Lt)$ and $s\geq j>i$.
\begin{algorithm}\mbox{}
\begin{enumerate}
    \item Remove the $a$-configuration of size $s-i$ from $T$ and
    slide it to get a $2\times i$ tableau.
    \item Remove the $1$'s, $\bar{1}$'s and the null-configuration 
    from the result to get a skew tableau of shape 
    $(i,i-t_{2T})/(t_{1T})$.
    \item Using the type $D$ sliding algorithm, produce a skew tableau
    of shape $((j),(j)-(t_{2T}+(j-i)))/(t_{1T}+(j-i))$.
    \item Fill this tableau with $1$'s, $\bar{1}$'s, and a 
    null-configuration so that the result is a $2\times j$ tableau.
    \item Use the height $2$ fill map $F_{2,s}$ (section \ref{sec:decomp}) to 
    insert $s-j$ columns into the tableau.
\end{enumerate}
\end{algorithm}
We can also define a map $\iota_{i}^{j}:B(i\Lt)\to B(j\Lt)\cup\{0\}$
for $j<i$ by
$$
\iota_{i}^{j}(T)=\begin{cases}
    (\iota_{j}^{i})^{-1}(T) & \textrm{ if }T\in\iota_{j}^{i}(B(j\Lt)), \\
    0 & \textrm{ otherwise }.
\end{cases}
$$
Reversing the above algorithm makes this map explicit.  Lastly, we
define $\iota_{i}^{i}$ to be the identity map on $B(i\Lt)$, so
$\iota_{i}^{j}$ is defined for all $i,j\in\{0,\ldots,s\}$.

We have already observed that by the $\dual$-duality of
$B(k\Lt)\subset\tbs$, each vertex $v_{\la}\in B(k\Lt)$ has a
complementary vertex $v'_{\la}\in B(k\Lt)$ such that
$R(v_{\la})+R(v'_{\la})=2s$.  We define the involution $\bcdual$ on
$\tbs$ as follows: Let $T\in B(v_{\la})$ such that
$T=\ftil{i_{1}}\cdots\ftil{i_{m}}u_{\la}$, where $u_{\la}$ is the
$U_q(D_{n-1})$-highest weight tableau of $B(v_{\la})$.  Then
$T^{\bcdual}=\ftil{i_{1}}\cdots\ftil{i_{m}}u'_{\la}$, where $u'_{\la}$
is the $U_q(D_{n-1})$-highest weight tableau of $B(v'_{\la})$.
Alternatively, this map is the composition of $\dual$ with the ``local
$\dual$'' map, which applies only to the tableaux in $B(v_{\la})$
viewed as a $U_q(D_{n-1})$-crystal.

We now define $\sigma(T)$ combinatorially.  Suppose $T\in
B(k\Lt)\subset\tbs$, and $\ell$ be minimal such that
$\iota_{k}^{s}(T)\in \iota_{\ell}^{s}(B(\ell\Lt))$.  Then
\begin{equation}\label{eq:sigma}
    \sigma(T)=\iota_{k}^{s+\ell-k}(T^{\bcdual})=(\iota_{k}^{s+\ell-k}(T))^{\bcdual},
\end{equation}
where it was used that $\iota_{i}^{j}$ commutes with $\bcdual$.

\subsection{Properties of $\ftil{0}$ and $\etil{0}$.}

This combinatorial approach immediately gives us useful information 
about this crystal, such as the following lemma.

\begin{Lemma} \label{lem:uk} For $k=0,1,\ldots,s$, let $u_{k}$ denote the highest
weight vector of the classical component $B(k\Lt)\subset\tbs$.  Then
$$
\ftil{0}(u_{k})=\begin{cases}
    u_{k+1} & \textrm{if }k<s, \\
    0 & \textrm{if }k=s.
\end{cases}
$$
\end{Lemma}

\begin{proof} Observe that
    $$
    u_{k}= \underbrace{
    \stack{1}{2}
    \cdots
    \stack{1}{2}}_{k}
    \underbrace{ 
    \stack{1}{\bar{1}}
    \cdots
    \stack{1}{\bar{1}}}_{s-k};
    $$
    We wish to calculate $\ftil{0}(u_{k})=\sigma\ftil{1}\sigma(u_{k})$.
    
    Note that $u_{k}\notin\iota_{k-1}^{k}(B((k-1)\Lt))$, so $\ell=k$ 
    in the combinatorial definition of $\sigma$ above.  It follows 
    that $\sigma(u_{k})=\iota_{k}^{s}(u_{k}^{\bcdual})$, which is
    \begin{equation*}
	\iota_{k}^{s}(u_{k}^{\bcdual}) 
	=
	\iota_{k}^{s}\Bigg( 
	\underbrace{
	\stack{1}{\bar{1}}
	\cdots
	\stack{1}{\bar{1}}}_{s-k}
	\underbrace{ 
	\stack{2}{\bar{1}}
	\cdots
	\stack{2}{\bar{1}}}_{k}\Bigg)
	=
	\framebox{$\emptyset_{s-k}$}
	\underbrace{ 
	\stack{2}{\bar{1}}
	\cdots
	\stack{2}{\bar{1}}}_{k},
    \end{equation*}
    where $\framebox{$\emptyset_{i}$}$ denotes a null-configuration of
    size $i$ (see Definition~\ref{def:null}).  
    If $k=s$, $\ftil{1}$ kills this tableau, as claimed in
    the second case of the Lemma.  Otherwise, acting by $\ftil{1}$
    will decrease the size of the null-configuration by $1$ and add
    another $\stack{2}{\bar{1}}$ to the columns on the right.  It
    follows that $\iota_{s}^{k}$ kills this tableau, but
    $\iota_{s}^{k+1}$ does not, so now $\ell=k+1$ in the combinatorial
    definition of $\sigma$.  Thus,
    \begin{equation*}
	\sigma\ftil{1}\sigma(u_{k})=\iota_{s}^{k+1}\Bigg(
	\underbrace{
	\stack{1}{2}
	\cdots
	\stack{1}{2}	}_{k+1}
	\framebox{$\emptyset_{s-k-1}$}
	\Bigg)
	=
	\underbrace{
	\stack{1}{2}
	\cdots
	\stack{1}{2}}_{k+1}
	\underbrace{ 
	\stack{1}{\bar{1}}
	\cdots
	\stack{1}{\bar{1}}}_{s-k-1}
	=u_{k+1}.
    \end{equation*}
\end{proof}

\begin{Cor} \label{cor:uk} Let $u_{k}$ be as above for $k>0$.  Then
    $$
    \etil{0}(u_{k})=	u_{k-1}     .
    $$
\end{Cor}

A similar combinatorial analysis can be carried out on lowest weight
tableaux to show that $\ftil{0}(u_{k}^{\dual})=u_{k-1}^{\dual}$ and
$\etil{0}(u_{k}^{\dual})=u_{k+1}^{\dual}$ for appropriate values of $k$.
Since $u_{0}=u_{0}^{\dual}$, this gives us the following Corollary:

\begin{Cor} \label{cor:phiuk} For highest weight vectors $u_{k}$ and lowest weight 
vectors $u_{k}^{\dual}$, we have
$$
\varphi_{0}(u_{k})=\vare_{0}(u_{k}^{\dual})=s-k
\qquad \text{and} \qquad
\varphi_{0}(u_{k}^{\dual})=\vare_{0}(u_{k})=s+k
$$
\end{Cor}

\section{Perfectness of $\tbs$}\label{sec:perfect}

\subsection{Overview} To show that $\tbs$ is perfect, it must be shown that
all criteria of Definition~\ref{def:perfect} are satisfied with $\ell=s$.  
We have taken part 3 of Definition~\ref{def:perfect} as part of our hypothesis 
for Theorem \ref{thm:main}, so we do not attempt to prove this here.

Part 2 of Definition~\ref{def:perfect} is satisfied by simply noting that
$\la=\om_2=s\Lambda_2-2s\Lz$ is a weight in $P_{\cl}$ such that $B_{\la}=\{u_{s}\}$
contains only one tableau and all other tableaux in $\tbs$ have
``lower'' weights.

In section \ref{sec:per con}, we show that $\tbs\otimes \tbs$ is
connected proving part 1 of Definition~\ref{def:perfect}.  Parts 4 and
5 of Definition~\ref{def:perfect} will be dealt with simultaneously in
sections~\ref{sec:sur} and~\ref{sec:inj} by examining the levels of
tableaux combinatorially.  We will see that the level of a generic
tableau is at least $s$ and the tableaux of level $s$ are in
bijection with the level $s$ weights.  In section~\ref{sec:unique} we
show that $\tbs$ is the unique affine crystal satsifying the
properties of Conjecture~\ref{conj:Brs} thereby proving
Theorem~\ref{thm:main}.

\subsection{Connectedness of $\tbs$} \label{sec:per con}

\begin{Lemma}[Part 1 of Definition~\ref{def:perfect}]
The crystal $\tbs\otimes \tbs$ is connected.
\end{Lemma}

\begin{proof} (This proof is very similar to that in~\cite[Proposition
5.1]{Koga:1999}.)  For $k=0,1,\ldots,s$, let $u_{k}$ denote the
highest weight vector of the classical component $B(k\Lt)\subset\tbs$,
as in Lemma \ref{lem:uk}.  We will show that an arbitrary vertex
$b\otimes b'\in\tbs\otimes \tbs$ is connected to $u_{0}\otimes u_{0}$.

We know that for some $j\in\{0,\ldots,s\}$, we have $b'\in B(j\Lt)$. 
Then for some pair of sequences $i_{1}, i_{2},\ldots,i_{p}$ (with
entries in $\{1,\ldots,n\}$) and $m_{1}, m_{2},\ldots,m_{p}$ (with
entries in $\mathbb{Z}_{>0}$) and some $b^{1}\in\tbs$, we have
$\etil{i_{1}}^{m_{1}}\etil{i_{2}}^{m_{2}}\cdots\etil{i_{p}}^{m_{p}}(b\otimes 
b')=b^{1}\otimes u_{j}$.

By Corollary \ref{cor:phiuk}, $\varphi_{0}(u_{j})=s-j$, so if
$\vare_{0}(b^{1})\leq s-j$, Lemma \ref{lem:uk} tells us that
$\etil{0}^{j}(b^{1}\otimes u_{j})=b^{1}\otimes u_{0}$.  If
$\vare_{0}(b^{1})=r>s-j$, then $\etil{0}^{r-s+j}(b^{1}\otimes
u_{j})=b^{2}\otimes u_{j}$, where $\vare_{0}(b^{2})=r-(r-s+j)=s-j$, so
$\etil{0}^{j}(b^{2}\otimes u_{j})=b^{2}\otimes u_{0}$.  In either
case, our arbitrary $b\otimes b'$ is connected to an element of the
form $b''\otimes u_{0}$.

Let $j'$ be such that $b''\in B(j'\Lt)$.  Since $u_{0}$ is the unique
element of $B(0)$, the crystal for the trivial representation of
$U_q(D_n)$, we know that $B(j'\Lt)\otimes B(0)\simeq B(j'\Lt)$. 
Therefore, $b''\otimes u_{0}$ is connected to $u_{j'}\otimes u_{0}$. 
Finally, we note that $\varphi_{0}(u_{0})=s<s+j'=\vare_{0}(u_{j'})$
for $j'\neq 0$, so $\etil{0}^{j'}(u_{j'}\otimes
u_{0})=\etil{0}^{j'}(u_{j'})\otimes u_{0}=u_{0}\otimes u_{0}$,
completing the proof.
\end{proof}

\subsection{Preliminary observations}\label{per15}
We first make a few observations. 

\begin{Prop} \label{prop:eps1 eps0}
    Let $T\in B(k\Lt)\subset\tbs$, and set $T_{m}=\iota_{k}^{m}(T)$
    for $m=s,s-1,\ldots,\ell$, where $\ell$ is minimal such that
    $\iota_{k}^{\ell}(T)\neq 0$.  If $\ell\neq s$, we have for 
    $s>m\geq\ell$
    \begin{eqnarray*}
	&&\vare_{1}(T_{m+1})=\vare_{1}(T_{m})+1 
	\textrm{ and } 
	\vare_{0}(T_{m+1})=\vare_{0}(T_{m})-1, \\
	&&\varphi_{1}(T_{m+1})=\varphi_{1}(T_{m})+1 
	\textrm{ and } 
	\varphi_{0}(T_{m+1})=\varphi_{0}(T_{m})-1.
    \end{eqnarray*}
\end{Prop}

\begin{proof}
    Let $\ell\leq m\leq s-1$, so $\iota_{m}^{m+1}$ is defined.  We
    first consider the difference between the reduced $1$-signatures
    of $D_{2,s}(T_{m})$ and $D_{2,s}(T_{m+1})=
    D_{2,s}(\iota_{m}^{m+1}(T_{m}))$, since the action of $\etil{1}$ on
    these tableaux is defined by the action of the classical
    $\etil{1}$ on their image under $D_{2,s}$.  Let $-^{M}+^{P}$ be the
    reduced $1$-signature of $D_{2,s}(T_{m})$.  As in section
    \ref{sec:edge}, let $r_{1}$ denote the number of $1$'s in
    $D_{2,s}(T_{m})$, $r_{3}$ the number of $\bar{1}$'s, $r_{2}$ the
    size of the null-configuration, and $t_{1}=r_{1}+r_{2}$,
    $t_{2}=r_{2}+r_{3}$.  Then there is a contribution $-^{r_{2}}+^{r_{2}}$
    to the 1-signature from the null-configuration, and the remaining $-$'s 
    and $+$'s come from $1$'s with a letter greater than $2$ below them and
    $\bar{1}$'s with a letter less than $\bar{2}$ above them,
    respectively.

    We now have two cases, as in section~\ref{sec:edge}.  If
    $t_{1}+t_{2}\geq s$, $\iota_{m}^{m+1}$ simply increases the size
    of the null-configuration in $D_{2,s}(T_{m})$ by $1$.  It follows
    that the reduced $1$-signature of $D_{2,s}(T_{m+1})$ is
    $-^{M+1}+^{P+1}$, as we wished to show.  On the other hand, if
    $t_{1}+t_{2}<s$, after step $(2)$ of Algorithm \ref{alg:iota} for
    $\iota_{m}^{m+1}$ we have a tableau of shape
    $(m,m-r_{3})/(r_{1})$.  In step (3), we slide this into shape
    $(m+1,m+1-(r_{3}+1))/(r_{1}+1)$.  We claim that the rightmost
    ``uncovered'' letter in the second row of this tableau is greater
    than $2$ and the leftmost ``unsupported'' letter in the first row
    is less than $\bar{2}$.  As observed in the preceding paragraph,
    this implies that after refilling the empty spaces as in step
    $(4)$ of Algorithm \ref{alg:iota} the reduced $1$-signature of our
    tableau is $-^{M+1}+^{P+1}$ in this case as well.
    
    Let us first consider the leftmost ``unsupported'' letter.  After 
    step $(2)$, our tableau is of the form
    $$
    \begin{array}{|c|c|c|c|c|c|c|c|c|}
        \cline{4-9}
         \multicolumn{2}{c}{} & & a_{r_{1}+1} & \cdots & a_{m-r_{3}} & 
         a_{m-r_{3}+1} & 
         \cdots & a_{s}  \\
        \cline{1-9}
        b_{1} & \cdots & b_{r_{1}} & b_{r_{1}+1} & \cdots & b_{m-r_{3}} &  
        \multicolumn{3}{c}{}  \\
        \cline{1-6}
    \end{array}
    $$
    and its column word is unchanged by the slide
    $$
    \begin{array}{|c|c|c|c|c|c|c|c|c|c|}
	\cline{4-5}	\cline{7-10}
	 \multicolumn{2}{c}{} & & a_{r_{1}+1} & \cdots & & a_{m-r_{3}} & 
	 a_{m-r_{3}+1} & 
	 \cdots & a_{s}  \\
	\cline{1-10}
	b_{1} & \cdots & b_{r_{1}} & b_{r_{1}+1} & \cdots & b_{m-r_{3}} &  
	\multicolumn{4}{c}{}  \\
	\cline{1-6}
    \end{array},
    $$    
    so we have $a_{m-r_{3}}<b_{m-r_{3}}\leq\bar{2}$.
    
    The second row of this tableau has $m-r_{3}$ boxes just as it did
    before sliding, so the boxes in the bottom row will never be
    moved.  It follows that this sliding procedure only changes
    L-shaped subtableaux into $\Le$-shapes \big(i.e., \quad $
    \begin{array}{|c|c|c}
        \cline{1-1}
         &  \multicolumn{2}{c}{}  \\
        \cline{1-2}
         &  &   \\
        \cline{1-2}
    \end{array}
    $
    into
    $
    \begin{array}{c|c|c|}
	\cline{3-3}
	   \multicolumn{1}{c}{} & & \\
	\cline{2-3}
	 &  &   \\
	\cline{2-3}
    \end{array}
    $
    \; \big) and never involves any $\Gamma$- or $\Lee$-shapes.
    According to the Lecouvey $D$-equivalence relations from
    section~\ref{sec:plac}, such moves can only be made when the
    letters in the bottom row are strictly greater than $2$.  
    Specifically, in relations $(3)$ and $(4)$, the letter which is 
    ``uncovered'' is either $n$ or $\bar{n}$, while in relations $(1)$ 
    and $(2)$ only the second case of each relation applies.  This 
    proves our claim, and thus the first half of the proposition.

Since $\etil{0}=\sigma\circ \etil{1}\circ\sigma$, we can derive the statements
about $\vare_0$ and $\varphi_0$ from the corresponding statements about
$\vare_1$ and $\varphi_1$. More precisely, $\vare_0(T)=\vare_1(\sigma(T))$
and $\varphi_0(T)=\varphi_1(\sigma(T))$ and by \eqref{eq:sigma} we have
\begin{equation*}
\sigma(T_m)=\bigl(\iota_m^{s+\ell-m}\circ\iota_k^m(T)\bigr)^{\bcdual}
=\iota_k^{s+\ell-m}(T^{\bcdual}).
\end{equation*}
Hence
\begin{multline*}
\vare_0(T_{m+1})=\vare_1(\sigma(T_{m+1}))
=\vare_1(\iota_k^{s+\ell-m-1}(T^{\bcdual}))\\
=\vare_1(\iota_k^{s+\ell-m}(T^{\bcdual}))-1
=\vare_1(\sigma(T_m))-1=\vare_0(T_m)-1.
\end{multline*}
A similar computation can be carried out for $\varphi_0$.
\end{proof}

\begin{Cor} \label{cor:lz plus lo}
    Given the above hypotheses, we have
    \begin{equation*}
	 \lan h_{0}+h_{1},\vare(T_{s})\ran= \lan
	 h_{0}+h_{1},\vare(T_{s-1})\ran= \cdots 
	 =\lan h_{0}+h_{1},\vare(T_{\ell})\ran\neq 0.
    \end{equation*}
\end{Cor}

The following observation is an immediate consequence of Remark 
\ref{rem:comb iota}:

\begin{Cor} \label{cor:hi}
    For $i=2,\ldots,n$,
    \begin{equation*}
	\lan h_{i},\vare(T_{s})\ran= \lan h_{i},\vare(T_{s-1})\ran= 
	\cdots  = \lan h_{i},\vare(T_{\ell})\ran.
    \end{equation*}
\end{Cor}

\begin{Lemma}\label{lem:ups}
    The map $\Upsilon_{s-1}^{s}:\mathcal{T}(s-1)\hookrightarrow\mathcal{T}(s)$
    viewed as sending the set underlying $\tilde{B}^{2,s-1}$ into the
    set underlying $\tbs$ increases the level of a tableau by exactly
    $1$.
\end{Lemma}

\begin{proof} This map is defined by sending each summand
$B(k\Lt)\subset\tilde{B}^{2,s-1}$ to $B(k\Lt)\subset\tbs$ for
$k=0,\ldots,s-1$, so
$\varphi_{i}(\Upsilon_{s-1}^{s}(T))=\varphi_{i}(T)$ for
$i=1,\ldots,n$.  To calculate the change in $\varphi_{0}(T)$, we must
consider the difference between $\varphi_{1}(\sigma_{s-1}(T))$ and
$\varphi_{1}(\sigma_{s}(\Upsilon_{s-1}^{s}(T)))$.  By our descriptions
of maps on crystals, we have
$\sigma_{s}(\Upsilon_{s-1}^{s}(T))=\Upsilon_{s-1}^{s}(\iota_{j}^{j+1}(\sigma_{s-1}(T)))$,
where $j$ is determined by $\sigma_{s-1}(T)\in
B(j\Lt)\subset\tilde{B}^{2,s-1}$.  By Proposition \ref{prop:eps1
eps0},
$\varphi_{1}(\iota_{j}^{j+1}(\sigma_{s-1}(T)))=\varphi_{1}(\sigma_{s-1}(T))+1$.
\end{proof}

\subsection{Surjectivity} \label{sec:sur}

Given a weight $\la\in(P_{\cl}^{+})_{s}$, we construct a tableau
$T_{\la}\in\tbs$ such that $\vare(T_{\la})=\varphi(T_{\la})=\la$.
This amounts to constructing $T_{\la}$ so that its reduced
$i$-signature is $-^{\vare_{i}(T_{\la})}+^{\vare_{i}(T_{\la})}$.  Note
that such a tableau is invariant under the $\dual$-involution, so its
symmetry allows us to define it beginning with the middle, and
proceeding outwards.

For $i=0,\ldots,n$, let $k_{i}=\lan h_{i},\la\ran$.  We first
construct a tableau $T_{\la'}$ corresponding to the weight
$\la'=\sum_{i=2}^{n}k_{i}\Lambda_{i}$.  We begin with the middle
$k_{n-1}+k_{n}$ columns of $T_{\la'}$.  If $k_{n-1}+k_{n}$ is even and
$k_n\geq k_{n-1}$, these columns of $T_{\la'}$ are
\begin{displaymath}
\underbrace{
\stack{n-2}{n-1}
\cdots
\stack{n-2}{n-1}
}_{k_{n-1}}
\underbrace{
\stack{n-1}{n}
\cdots
\stack{n-1}{n}
}_{(k_{n}-k_{n-1})/2}
\underbrace{
\stack{\overline{n}}{\overline{n-1}}
\cdots
\stack{\overline{n}}{\overline{n-1}}
}_{(k_{n}-k_{n-1})/2}
\underbrace{
\stack{\overline{n-1}}{\overline{n-2}}
\cdots
\stack{\overline{n-1}}{\overline{n-2}}
}_{k_{n-1}}
\end{displaymath}
If $k_{n-1}+k_{n}$ is odd and $k_{n}\geq k_{n-1}$, we have
\begin{displaymath}
\underbrace{
\stack{n-2}{n-1}
\cdots
\stack{n-2}{n-1}
}_{k_{n-1}}
\underbrace{
\stack{n-1}{n}
\cdots
\stack{n-1}{n}
}_{(k_{n}-k_{n-1}-1)/2}
\stack{\bar{n}}{n}
\underbrace{
\stack{\overline{n}}{\overline{n-1}}
\cdots
\stack{\overline{n}}{\overline{n-1}}
}_{(k_{n}-k_{n-1}-1)/2}
\underbrace{
\stack{\overline{n-1}}{\overline{n-2}}
\cdots
\stack{\overline{n-1}}{\overline{n-2}}
}_{k_{n-1}}
\end{displaymath}
In either case, if $k_{n}<k_{n-1}$, interchange $n$ and
$\bar{n}$, and $k_n$ and $k_{n-1}$ in the above configurations.

Next we put a configuration of the form
$$
\underbrace{
\stack{1}{2}
\cdots
\stack{1}{2}
}_{k_{2}}
\underbrace{
\stack{2}{3}
\cdots
\stack{2}{3}
}_{k_{3}}
\qquad
\cdots
\qquad
\underbrace{
\stack{n-3}{n-2}
\cdots
\stack{n-3}{n-2}
}_{k_{n-2}}
$$
on the left, and a configuration of the form
$$
\underbrace{
\stack{\overline{n-2}}{\overline{n-3}}
\cdots
\stack{\overline{n-2}}{\overline{n-3}}
}_{k_{n-2}}
\underbrace{
\stack{\overline{n-3}}{\overline{n-4}}
\cdots
\stack{\overline{n-3}}{\overline{n-4}}
}_{k_{n-3}}
\qquad
\cdots
\qquad
\underbrace{
\stack{\bar{2}}{\overline{1}}
\cdots
\stack{\bar{2}}{\overline{1}}
}_{k_{2}}
$$
on the right.  Denote the set of tableaux constructed by the procedure
up to this point by $\T(s')$.

Observe that the reduced $1$-signature of $T_{\la'}$ is empty, so
$\lan h_{1},\varphi(T_{\la'})\ran=0$.  Furthermore, since $\la'$ has
the same number of $1$'s as $\bar{1}$'s, it is fixed by $\sigma$, so
$\lan h_{0},\varphi(T_{\la'})\ran=0$ as well.  Thus Proposition
\ref{prop:eps1 eps0} implies that $T_{\la'}\in \tbs_{\min}\cap
B(s'\Lt)\setminus\iota_{s'-1}^{s'}(B((s'-1)\Lt))$ as a subset of
$\tilde{B}^{2,s'}$.  Recall the embedding
$\Upsilon_{s'}^{s}:\tilde{B}^{2,s'} \hookrightarrow\tbs$ from
Definition~\ref{def:upsilon}.  Since $s'$ is the minimal $m$ for which
$\iota_{s'}^{m}(T_{\la'})\neq 0$, Lemma \ref{lem:ups} and its proof
tell us that $\vare_{0}(F_{2,s}(T_{\la'}))=s-s'$, where the fill map
$F_{2,s}$ inserts an $a$-configuration to increase the width of
$T_{\la'}$ to $s$.  By the same proposition the desired tableau is
$T_{\lambda}=\iota_{s'}^{s'+k_1}\circ F_{2,s}(T_{\lambda'})$.  We
denote by $\Mmin{s}$ the set of tableaux constructed by this
procedure.

\subsection{Injectivity}\label{sec:inj}

In this subsection we show that the tableaux in $\Mmin{s}$ are all the
minimal tableaux in $\tbs$.

We first introduce some useful notation.  Observe that any tableau $T$
can be written as $T=T_{1}T_{2}T_{3}T_{4}T_{5}$, where the block
$T_{i}$ has width $k_{i}$, and all letters in $T_{1}$ (resp.  $T_{5}$)
are unbarred or $\bar{n}$ in the second row (resp.  barred or $n$ in
the first row), all columns in $T_{2}$ (resp.  $T_{4}$) are of the
form $\stack{a}{\bar{b}}$ where $a<b\leq n-1$ (resp.  $b<a\leq n-1$),
and all columns in $T_{3}$ are of the form $\stack{a}{\bar{a}}$ for
some $a$.  Note that for a tableau in $B(s\Lt)\subset\tbs$ we have
$0\leq k_{3}\leq 1$.  Also note that $T_{2}$ and $T_{4}$ do not contain 
any $n$'s or $\bar{n}$'s.

\begin{theorem} \label{thm:injective}
    We have $\Mmin{s}=\tbs_{\min}$.
\end{theorem}

\begin{proof}
    We prove the theorem by induction on $s$.  For the base case, we
    checked explicitly that the statement of the theorem is true for
    $s=0,1,2$. 
    
    By our induction hypothesis,
    $\Mmin{s-1}=\tilde{B}^{2,s-1}_{\min}$.  By Lemma \ref{lem:ups},
    $\Upsilon_{s-1}^{s}$ increases the level of a tableau by $1$, and
    therefore
    $\Upsilon_{s-1}^{s}(\tilde{B}^{2,s-1}_{\min})=(\tbs\setminus
    B(s\Lt))_{\min}$.  By Corollaries \ref{cor:lz plus lo} and
    \ref{cor:hi}, $\iota_{s-1}^{s}$ does not change the level of a
    tableau, so it suffices to show that
    \begin{equation}     \label{eq:minimal}
    \T(s)=\tbs_{\min}\cap(B(s\Lt)\setminus\iota_{s-1}^{s}(B((s-1)\Lt))).
    \end{equation}
    By Lemmas \ref{lem:plus minus} and \ref{lem:mixed} below, if $T$ is a
    minimal tableau not in the image of $\iota_{s-1}^{s}$, it has
    $k_{2}=k_{4}=0$, and if $k_{3}=1$, then $T_{3}=\stack{n}{\bar{n}}$
    or $T_{3}=\stack{\bar{n}}{n}$.  By Lemma \ref{lem:unmixed},
    equation \eqref{eq:minimal} follows.
\end{proof}

Here we state the lemmas used in the proof of Theorem~\ref{thm:injective}.
The proofs are given in Appendices~\ref{app:A}-\ref{app:D} and together
with Theorem~\ref{thm:injective} all rely on induction on $s$.
The base cases $s=0,1,2$ have been checked explicitly.

\begin{Lemma}\label{lem:lb}
For all $t\in \tbs$ we have $\lan c,\varepsilon(t)\ran \ge s$ and 
$\lan c,\varphi(t)\ran \ge s$.
\end{Lemma}

\begin{Lemma} \label{lem:plus minus}
If $T\in (\tbs\cap B(s\Lt))_{\min}$, we have $k_{1}+k_{2}=k_{4}+k_{5}=\lfloor
s/2 \rfloor$ with $k_1,k_2,k_4$ and $k_5$ as defined in the beginning of this
subsection.
\end{Lemma}

\begin{Lemma} \label{lem:mixed}
    Suppose $T\in (\tbs\cap B(s\Lt))_{\min}$ has both an unbarred
    letter and a barred letter in a single column $\stack{a}{
    \bar{b}}$ other than $\stack{n}{\bar{n}}$, $\stack{\bar{n}}{ n}$,
    $\stack{n-1}{\bar{n}}$, or $\stack{n}{ \overline{n-1}}$.  Then
    $T\in\iota_{s-1}^{s}(B((s-1)\Lt)$.
\end{Lemma}

\begin{Lemma} \label{lem:unmixed}
Let $T\in (\tbs\cap B(s\Lt))_{\min}$ such that $T$ does not contain
any column $\stack{a}{\bar{b}}$ for $1\le a,b\le n$ except possibly
$\stack{n-1}{\bar{n}}$, $\stack{n}{\bar{n}}$, or
$\stack{n}{\overline{n-1}}$.  Then $T\in \T(s)$.
\end{Lemma}

\subsection{Uniqueness} \label{sec:unique}
Theorem~\ref{thm:main} follows as a corollary from the next
Proposition.

\begin{Prop}\label{prop:uniqueness}
$\tbs$ is the unique affine finite-dimensional crystal structure satisfying
the properties of Conjecture~\ref{conj:Brs}.
\end{Prop}

For the proof of Proposition~\ref{prop:uniqueness} we must show that
our choice of $\sigma$ is the only choice satisfying the properties of
Conjecture~\ref{conj:Brs}.  Recall from the beginning of
section~\ref{sec:aff} the relationship between $\sigma$ and $\hsig$.
Let $T\in\tbs$.  We know that
\begin{equation}
    \label{eq:sigma prop}
    \wt(T)=\sum_{i=0}^{n}k_{i}\Lambda_{i}\Leftrightarrow 
    \wt(\sigma(T))=k_{1}\Lz+k_{0}\Lo+\sum_{i=2}^{n}k_{i}\Lambda_{i}.
\end{equation}
Once $\hsig$ is determined, the given definition of $\sigma$ (sending
$D_{n-1}$ highest weight vectors to $D_{n-1}$ highest weight vectors,
etc.)  is the only involution of the set of tableaux in $\tbs$
satisfying \eqref{eq:sigma prop} and agreeing with $\hsig$.

As we observed in section~\ref{sec:hsig}, $\hsig(v)$ and $v$ must be
associated with the same partition, and if $v'$ is the complementary
vertex of $v$, $\hsig(v')$ is the complementary vertex of $\hsig(v)$.
We now prove a few lemmas that uniquely determine $\hsig$.

Please note that in this section we often use the phrase ``the tableau
$b$ is in the branching component vertex $v$'' to mean $b\in
B(v)$.

\begin{Lemma}
    Suppose $b\in\tbs$ is in a branching component vertex of rank $p$
    and $\ftil{0}(b)\neq 0$.  Then the branching component vertex
    containing $\ftil{0}(b)$ has rank $p-1$.
\end{Lemma}

\begin{proof}
    Recall from the weight structure of type $D$ algebras that
    $\alpha_{0}=2\Lz-\Lambda_2$ and $\wt(\ftil{0}(b))=\wt(b)-\alpha_{0}$.
    Define $\cw(b)=\wt(b)-(\varphi_{0}(b)-\vare_{0}(b))\Lambda_0$.
    The above implies that
    $\cw(\ftil{0}(b))=\cw(b)+\Lambda_2=\cw(b)+\varepsilon_1+\varepsilon_2$.
    Similarly, $\cw(\ftil{i}(b))=\cw(b)-\alpha_i$, so that by
    \eqref{eq:alpha} only $\ftil{1}$ changes the $\varepsilon_1$
    component by $-1$.  Since $\ftil{1}$ increases the rank by one and
    $\ftil{0}$ changes the $\varepsilon_1$ component by $+1$, it
    follows that $\ftil{0}$ decreases the rank by one.
\end{proof}

Since $\ftil{1}$ increases the rank by one, $\ftil{0}$ decreases the rank
by one and $\ftil{0}=\sigma\ftil{1}\sigma$, the following corollary holds.

\begin{Cor}\label{cor:flip}
    Suppose $b$ is in a branching component vertex $v$ of rank $p$.
    Then $\hsig(v)$, the branching component vertex containing
    $\sigma(b)$, has rank $2s-p$.
\end{Cor}

Note that this determines $\hsig$ on
$B(s\Lt)\setminus\iota_{s-1}^{s}(B((s-1)\Lt))$, and is in agreement
with our definition of $\hsig$ restricted to that domain.

\begin{Lemma} \label{lem:middle}
    Let $v\in\Sup(k\Lt)$ be a branching component vertex of rank $s$
    associated with a rectangular partition, and let $\ell$ be minimal such that
    $v\in\iotil_{\ell}^{k}(\Sup(\ell\Lt))$.  Then the
    hypothesis that $\tbs$ is perfect of level $s$ can only be
    satisfied if $\hsig(v)$ is the vertex associated with the same
    shape as $v$ with rank $s$ in $\Sup((s+\ell-k)\Lt)$.
\end{Lemma}

\begin{proof} 
    We have already shown that $\hsig(v)$ has the same shape as $v$
    and has rank $s$, so it only remains to show that
    $\hsig(v)\in\Sup((s+\ell-k)\Lt)$.
    
    First, observe that $v$ must contain a minimal tableau as
    constructed in section \ref{sec:sur}, according to the following 
    table.
    
    \vspace{.1in}    
    \begin{center}
    $
    \begin{array}{c|c}
	\textrm{shape associated with }v & \textrm{weight of tableau
	in }v \\
	\hline
	(2m) & m\Lambda_2+(k-2m)\Lo+(s-k)\Lz \\
	(m,m)  & m\Lambda_{3}+(k-2m)\Lo+(s-k)\Lz \\
    \end{array}
    $
    \end{center}
    \vspace{.1in}

    Let $T$ be the tableau constructed by this prescription, so that
    $\lan c,\cw(T)\ran=k$.  The criterion that $\lan
    c,\varphi(T)\ran\geq s$ forces us to have
    $\varphi_{0}(T)=\varphi_{1}(\sigma(T))\geq s-k.$ We denote
    $T_{i}=\iota_{k}^{i}(T)$, and thus have
    $\varphi_{1}(T_{i})=i-\ell$.  
    
    We show inductively that $\sigma(T_{i})=T_{s+\ell-i}$ for
    $\ell\leq i\leq s$.  As a base case, we see that $\lan
    c,\cw(T_{\ell})\ran=\ell$, so we must have
    $\varphi_{1}(\sigma(T_{\ell}))\geq s-\ell$.  The only $T_{i}$ for
    which this inequality holds is $T_{s}$, where we have
    $\varphi_{1}(T_{s})=s-\ell$.
    
    For the induction step, assume that $\sigma$ sends
    $T_{\ell},T_{\ell+1},\ldots,T_{k-1}$ to
    $T_{s},\ldots,T_{s+\ell-k+1}$, respectively.  By the above
    inequality this implies $\varphi_{1}(\sigma(T_{k}))\geq s-k$,
    which specifies that $\sigma(T_{k})=T_{s+\ell-k}$.
\end{proof}

\begin{Def}\label{def:ups check}
    Recall the definition of $\Upsilon_{s'}^{s}$ from
    Definition~\ref{def:upsilon}.  We define
    $\check{\Upsilon}_{s'}^{s}:\Sup(\tilde{B}^{2,s'})\hookrightarrow
    \Sup(\tilde{B}^{2,s})$ by $\check{\Upsilon}_{s'}^{s}(v)=v'$ if for
    some $T\in B(v)$, we have $\Upsilon_{s'}^{s}(T)\in B(v')$.
\end{Def}

\begin{Lemma} \label{lem:down}
    Let $v\in\Sup(k\Lt)$ be a branching component vertex of rank
    $1\leq p\leq s-1$ associated to the shape
    $(\lambda_{1},\lambda_{2})$. Suppose that for the branching
    component vertex $w\in\Sup(k\Lt)$ of rank $p+1$ with shape
    $(\lambda_{1}-1,\lambda_{2})$, $\tbs$ has the correct energy
    function and is perfect only if $\hsig(w)$ is as described in
    section \ref{sec:hsig}.  Then $\tbs$ has the correct energy
    function only if $\hsig(v)$ is as described in section
    \ref{sec:hsig}.
\end{Lemma}

\begin{proof}
    First, recall that the partitions associated to vertices of rank
    $p$ in $\Sup(\tbs)$ are produced by adding or removing one box
    from the partitions associated to vertices of rank $p-1$.  Since
    the vertex of rank $0$ is associated with a rectangle of shape $(s)$,
    the lowest rank for which we can have a two-row rectangle is $s$.
    It follows that removing a box from the first row of $v$ results
    in a partition of rank $p+1$, so there is in fact a vertex
    $w$ as described in the statement of the lemma.
        
    Let $\ell$ be minimal such that
    $v\in\iotil_{\ell}^{k}(\Sup(\ell\Lt))$.  We may assume
    $v\notin\Sup(\ell\Lt)$, and let $\hsig(v)$ be determined by the
    involutive property of $\hsig$ in the case $v\in\Sup(\ell\Lt)$.
    Specifically, we will show that the vertex $v'\in\Sup(s\Lt)$ with
    the same shape and rank as $v$ has the property that $\hsig(v')$
    is the complementary vertex of $v$, and therefore $\hsig(v)$ is
    the complementary vertex of $v'$.
    
    The top row of $w$ is one box shorter than the top row of
    $v$, so that $w\in\iotil_{\ell-1}^{k}(\Sup((\ell-1)\Lt))$, and
    $\ell-1$ is minimal with this property.  By our hypothesis,
    $\hsig(w)$ is the vertex with shape $(\lambda_{1}-1,\lambda_{2})$
    of rank $2s-p-1$ in $\Sup((s+(\ell-1)-k)\Lt)$.
    
    We now use induction on $s$.  Suppose that the only choice of
    $\hsig$, for which $\tilde{B}^{2,s-1}$ is perfect and has an energy
    function, is the choice of section~\ref{sec:hsig}.  Part 3 of
    Conjecture~\ref{conj:Brs} states that in $\tbs$, the energy on the
    component $B((s-k)\Lt)$ is $-k$, and so the difference in energy
    between $B((s-k)\Lt)$ and $B((s-j)\Lt)$ is $j-k$.  In order for
    this to be true for all $1\leq k,j\leq s-1$, the action of
    $\ftil{0}$ and $\etil{0}$ on $\tbs$ must agree with the action on
    $\tilde{B}^{2,s-1}$.  More precisely, if $v$ and $w$ are in
    different classical components of $\tilde{B}^{2,s-1}$ and
    $\ftil{0}(v)=w$ in $\tilde{B}^{2,s-1}$, then
    $\ftil{0}(\Upsilon_{s-1}^{s}(v))=\Upsilon_{s-1}^{s}(w)$; this
    statement extends naturally to $\Sup(\tilde{B}^{2,s-1})$ and
    $\Sup(\tbs)$.
    
    Let $v^{\dagger}$ denote the vertex with shape
    $(\lambda_{1},\lambda_{2})$ of rank $2s-p$ in
    $\Sup((s+\ell-k)\Lt)$.  Since we assumed $k\neq\ell$, we know that
    $v^{\dagger}\notin\Sup(s\Lt)$, and therefore $v^{\dagger}$ has a
    preimage under $\check{\Upsilon}_{s-1}^{s}$.  From our
    construction of $\hsig$ we know that in $\Sup(\tilde{B}^{2,s-1})$,
    $(\check{\Upsilon}_{s-1}^{s})^{-1}(v^{\dagger})$ has a $0$ arrow
    to $(\check{\Upsilon}_{s-1}^{s})^{-1}(\hsig(w))$.  Our induction
    argument tells us that in $\Sup(\tbs)$, $v^{\dagger}$ has a $0$
    arrow to $\hsig(w)$.  Since $v$ has a $1$ arrow to $w$ and we must
    have $\ftil{0}=\sigma\ftil{1}\sigma$, we conclude that in fact
    $\hsig(v)=v^{\dagger}$.
\end{proof}

\begin{Lemma} \label{lem:up}
    Let $v\in\Sup(k\Lt)$ be a branching component vertex of rank $s$
    associated to the non-rectangular shape
    $(\lambda_{1},\lambda_{2})$, and suppose that for the branching
    component vertex $w\in\Sup(k\Lt)$ of rank $s-1$ with shape
    $(\lambda_{1},\lambda_{2}+1)$, $\tbs$ has the correct energy
    function and is perfect only if $\hsig(w)$ is as described in
    section \ref{sec:hsig}.  Then $\tbs$ has the correct energy
    function only if $\hsig(v)$ is as described in section
    \ref{sec:hsig}.
\end{Lemma}

\begin{proof}
    (This proof is very similar to the proof of Lemma~\ref{lem:down}.)
    
    Let $\ell$ be minimal such that
    $v\in\iotil_{\ell}^{k}(\Sup(\ell\Lt))$, assuming $k\neq\ell$.
    Note that $\ell$ is also minimal for
    $w\in\iotil_{\ell}^{k}(\Sup(\ell\Lt))$, since the shapes for $v$
    and $w$ have the same number of boxes in the first row.  By our
    hypothesis, $\hsig(w)$ is the vertex with shape
    $(\lambda_{1},\lambda_{2}+1)$ of rank $s+1$ in
    $\Sup((s+\ell-k)\Lt)$.  Let $v^{\dagger}$ be the vertex with shape
    $(\lambda_{1},\lambda_{2})$ of rank $s$ in $\Sup((s+\ell-k)\Lt)$.
    From our construction of $\hsig$, we know that in
    $\Sup(\tilde{B}^{2,s-1})$, $(\check{\Upsilon}_{s-1}^{s})^{-1}(w)$
    has a $0$ arrow to
    $(\check{\Upsilon}_{s-1}^{s})^{-1}(\hsig(v^{\dagger}))$.  It
    follows from our induction argument that in $\Sup(\tbs)$, $w$ has
    a $0$ arrow to $\hsig(v^{\dagger})$.  Since $w$ has a $1$ arrow to
    $v$ and we must have $\ftil{0}=\sigma\ftil{1}\sigma$, we conclude
    that in fact $\hsig(v)=v^{\dagger}$.
\end{proof}

\begin{Cor} Corollary~\ref{cor:flip} and Lemmas~\ref{lem:middle},
\ref{lem:down} and \ref{lem:up} determine $\hsig$ on $\Sup(\tbs)$
uniquely.
\end{Cor}

\begin{proof}
For any vertex $v$ associated with shape
$(\lambda_{1},\lambda_{2})$ with rank $p\leq s$, 
$\hsig(v)$ is fixed by the image under $\hsig$ of a vertex with shape
$(\lambda_{1}-s+p,\lambda_{2})$ and rank $s$ by Lemma~\ref{lem:down}.  If
$\lambda_{1}-s+p\neq\lambda_{2}$, Lemmas~\ref{lem:down}
and~\ref{lem:up} may be used together to reduce determining $\hsig(v)$
to determining the action of $\hsig$ on a rectangular vertex of rank
$s$, which is given by Lemma~\ref{lem:middle}.
\end{proof}

\section{Discussion}\label{sec:discussion}

In this section we discuss some applications and open problems
regarding the crystals $\tbs$ introduced in this paper.

The major open question regarding $\tbs$ is of course its existence,
which was assumed throughout this paper.  A possible method of proof
is a generalization of the fusion construction of~\cite{KKMMNN:1992}.

In \cite{Kash:2003}, Kashiwara conjectures that for any quantum affine algebra,
$B^{r,s}$ is isomorphic as a classical crystal to a Demazure
crystal in an irreducible affine highest weight module of weight
$s\max(1,2/(\alpha_r,\alpha_r))\varpi_r-s\Lambda_0$, where
$\varpi_r=\Lambda_r-a^\vee_r\Lambda_0$ (except for type $A_{2n}^{(2)}$).
The $\ftil{0}$ edges that stay within the Demazure crystal should be among the
$\ftil{0}$ edges of $B^{r,s}$. The combinatorial structure of $\tbs$ as
constructed in this paper might give a hint on how to make this correspondence
more precise.

For a tensor product of affine finite crystals $B=B_L\otimes\dotsm\otimes B_1$ and 
a dominant integral weight $\la$ define the set of classically restricted paths as
\begin{equation*}
\Path(B,\la)=\bigl\{ b\in B \mid 
 \text{$\wt(b)=\la$, $\etil{i}(b)=0$ for all $i\in J$} \bigr\}
\end{equation*}
where $J=\{1,2,\ldots,n\}$.
The classically restricted one dimensional sum is defined to be
\begin{equation*}
X(B,\la;q)=\sum_{b\in\Path(B,\la)} q^{D_B(b)}.
\end{equation*}

In~\cite[Section 4]{HKOTT:2001} fermionic formulas $M(B,\la;q)$ are
defined which are sums of products of $q$-binomial coefficients, and it
is conjectured that $X(B,\la;q)=M(B,\la;q)$. This conjecture has been
proven for type $A_n^{(1)}$ \cite{KKR:1988,KR:1988,KSS:2002} and various other 
cases~\cite{Chari,Na:2001,OSS:2002,OSS:2003,OSS:2003a,S:2004,SS:2004}.
It is expected that the $X=M$ conjecture can also be proven in the 
case of tensor products of crystals $\tbs$ as constructed in this paper
by using the splitting map (see \cite{SS:2004}) and the single column
bijection of type $D_n^{(1)}$ (see \cite{S:2004}).  Using these maps, one should obtain a statistic-preserving bijection between the set of tableaux $\mathcal{T}(s)$ defined in section \ref{sec:decomp} and a set of rigged configurations which naturally indexes the $q$-binomial coefficients in the fermionic formulas.  The statistics that are preserved by this bijection (energy in the case of tableaux, co-charge in the case of rigged configurations) are the exponents of $q$ in the $X=M$ formula.

\appendix

\section{Proof of Lemma~\ref{lem:lb}} \label{app:A}
By induction hypothesis, $\tilde{B}^{2,s-1}_{\min}=\Mmin{s-1}$.

Observe that by Corollaries \ref{cor:lz plus lo} and \ref{cor:hi}
$\iota_{i}^{j}$ is level preserving, and by Lemma \ref{lem:ups}, the
map $\Upsilon_{s-1}^{s}$ increases the level of a tableau by one.  Our induction
hypothesis therefore allows us to assume that $t\in B(s\Lt)\setminus
\iota_{s-1}^{s}(B((s-1)\Lt)$.  Combinatorially, we may characterize
such tableaux as being those which are legal in the classical sense
and for which removing all $1$'s, $\bar{1}$'s, and null configurations
produces a tableau which is Lecouvey $D$-equivalent to a tableau whose
first row has width $s$.  This characterization follows from the 
combinatorial description of $\iota_{s-1}^{s}$ in Algorithm~\ref{alg:iota}.

We may further restrict our attention by the observation that if $T$ 
is minimal, then so is $T^{\dual}$.  We may therefore assume $T$ to be 
in the top half (inclusive of the middle row) of the branching 
component graph.  This means that $T$ has no more $\bar{1}$'s than $1$'s.

Our approach is to consider the tableau $T'$ that
results from removing the leftmost column from $T$.  We will show that
if $T'$ is minimal, the level of $T$ exceeds the level of $T'$ by at
least $2$, and if $T'$ is not minimal, the level of $T$ is at least as
great as the level of $T'$.

First consider the case when $T'$ is minimal.  Since $T$ is assumed to
be such that removing all $1$'s and $\bar{1}$'s produces a tableau
which is Lecouvey $D$-equivalent to a tableau whose first row has
width $s$, it is the case that removing all $1$'s and $\bar{1}$'s from
$T'$ produces a tableau which is Lecouvey $D$-equivalent to a tableau
whose first row has width $s-1$.  The minimal tableaux of
$\tilde{B}^{2,s-1}$ with this property are precisely $\T(s-1)$.  By
properties of $\T(s-1)$, we know that $\varphi_{0}(T')=0$, so
$\varphi_{0}(T)\geq\varphi_{0}(T')$.  Since our base case is $s=2$, we
know that the first column of $T$ is $\stack{a}{b}$, where $a$ and $b$ 
are both unbarred.  Observe that 
$\varphi(T)=\varphi(T')+2\Lambda_{b}+\textrm{non-negative weight}$ if
$b\neq n-1$ or  
$\varphi(T)=\varphi(T')+\Lambda_{n-1}+\Lambda_{n}+\textrm{non-negative 
weight}$ if $b=n-1$.  Hence, the level is increased by at least $2$.

Now suppose $T'$ is not minimal.  The level of the $i$-signatures
(that is to say, the level of the sum of the weights $\varphi_{i}(T')$
which depend on $i$-signatures) cannot have a net decrease for
$i=1,\ldots,n$, but there is now a possibility that
$\varphi_{0}(T)<\varphi_{0}(T')$.  We will show that when
$\varphi_{0}(T)<\varphi_{0}(T')$, the level of the $i$-signatures goes
up by at least $\varphi_{0}(T')-\varphi_{0}(T)$.

First, suppose $T$ has no $1$'s.  Then by one of our hypotheses, it
also has no $\bar{1}$'s, and is therefore fixed by $\sigma$: it
follows that $\varphi_{0}(T)=\varphi_{1}(T)$, so we may assume the 
upper-left entry of $T$ to be $1$.

We know that $\varphi_{0}(T)$ is equal to the number of $-$'s in the 
reduced $1$-signature of $\sigma(T)$.  Consider the following 
tableaux:
\begin{eqnarray*}
    T' & = & \underbrace{\stack{1}{b_{1}}\cdots\stack{1}{
    b_{m_{1}}}}_{m_{1}}
    \stack{a_{m_{1}+1} }{ b_{m_{1}+1}}\cdots\stack{a_{s-m_{2}} }{ b_{s-m_{2}}}
    \underbrace{\stack{a_{s-m_{2}+1} }{ \bar{1}}\cdots\stack{a_{s} }{ 
    \bar{1}}}_{m_{2}}\\
    T & = &  \stack{1 }{b}\underbrace{\stack{1}{b_{1}}\cdots\stack{1 }{ 
    b_{m_{1}}}}_{m_{1}}
    \stack{a_{m_{1}+1} }{ b_{m_{1}+1}}\cdots\stack{a_{s-m_{2}} }{ b_{s-m_{2}}}
    \underbrace{\stack{a_{s-m_{2}+1} }{\bar{1}}\cdots\stack{a_{s}}{ 
    \bar{1}}}_{m_{2}}\\
    \sigma(T') & = & \underbrace{\stack{1 }{ 
    b_{1}}\cdots\stack{1}{b_{m_{2}-1}}\stack{1 }{ 
    b'_{m_{2}}}}_{m_{2}}
    \stack{a'_{m_{2}+1} }{ b'_{m_{2}+1}}\cdots\stack{a'_{s-m_{1}} }{ 
    b'_{s-m_{1}}}
    \underbrace{\stack{a'_{s-m_{1}+1} }{ \bar{1}}\cdots\stack{a'_{s} }{ 
    \bar{1}}}_{m_{1}} \\
    \sigma(T) & = & \underbrace{\stack{1 }{ b}\stack{1 }{ b_{1}}\cdots\stack{1 }{ 
    b_{m_{2}-1}}}_{m_{2}}
    \stack{a''_{m_{2}} }{ b''_{m_{2}}}\cdots\stack{a''_{s-m_{1}-1} }{ 
    b''_{s-m_{1}-1}}
    \underbrace{\stack{a''_{s-m_{1}} }{ \bar{1}}\cdots\stack{a''_{s}}{ 
    \bar{1}}}_{m_{1}+1} .
\end{eqnarray*}
Note that our assumption that $m_{2}\leq m_{1}+1$ ensures that the 
absence of primes on $b,b_{1},\ldots,b_{m_{2}-1}$ is accurate.  

Let us consider all possible ways for the number of $-$'s in the 
$1$-signature to be smaller for $\sigma(T)$ than for $\sigma(T')$.  The 
number of $1$'s is the same, so the only way this contribution could 
be decreased is by having more $2$'s in the first $m_{2}$ letters of 
the bottom row.  This can only come about by having $b=2$, and only 
one $-$ may be removed in this way.

The other possibility is for the number of $-$'s contributed by
$\bar{2}$'s to be decreased.  The only Lecouvey relation which removes
a $\bar{2}$ assumes the presence of a column $\stack{2}{\bar{2}}$,
which we disallow (null-configuration).  To decrease this contribution
therefore requires an additional $+$ in the $1$-signature of
$\sigma(T)$ compared to that of $\sigma(T')$, which will bracket one
of the $-$'s from a $\bar{2}$.  The additional $+$ may come from one
of the additional $\bar{1}$'s, or from a $2$ that is ``pushed out''
from under the $1$'s at the beginning in the case $b=2$.  Note that
this second possibility is mutually exclusive with having more $2$'s
bracketing $1$'s at the beginning.

In any case, we see that $\varphi_{0}(T)-\varphi_{0}(T')\leq2$, and
that when this value is $2$, the first column of $T$ is $\stack{1}{ 2}$.
This column adds no $+$'s to the $i$-signatures, but does provide a
new $-$ in the $2$-signature.  Since $\Lambda_2$ is a level $2$ weight, the
level stays the same in this case.

If $\varphi_{0}(T)-\varphi_{0}(T')=1$ and the first column of $T$ is 
$\stack{1}{ 2}$, we in fact have a net increase in level.  If $b\neq 2$, 
the $i$-signature levels go up by at least 1, so still the total level 
cannot decrease.

\section{Proof of Lemma~\ref{lem:plus minus}} \label{app:B}
We first establish that $\lan c,\varphi(T)\ran\geq
2k_{1}+2k_{2}+k_{3}$, and thus by $\dual$-duality, $\lan
c,\vare(T)\ran\geq k_{3}+2k_{4}+2k_{5}$ as well.  Recall that $0\leq 
k_{3}\leq 1$.

First, observe that every letter in the bottom row of $T_{1}$
contributes:
\begin{itemize}
    \item a $-$ to the
    reduced $a$-signature if $2\leq a\leq n-2$ is in the bottom row;

    \item a $-$ to both the $(n-1)$-signature and the $n$-signature if
    $n-1$ is in the bottom row;

    \item a $-$ to the $n$-signature (resp.  $(n-1)$-signature) if $n$
    (resp.  $\bar{n}$) is the bottom row.
\end{itemize}
Suppose $T_{1}$ has a column of the form $\stack{a}{b}$ with $b\neq
a+1$, or $b=\bar{n}$ and $a\neq n-1$.  For the $-$ in the
$a$-signature of $T$ contributed by this $a$ to be bracketed, we must
have a column of the form $\stack{a'}{a+1}$ to the left of this column
in $T_{1}$, with $a'<a$.  Applying this observation recursively, we
see that to bracket as many $-$'s as possible we must eventually have
a column of the form $\stack{1}{c}$ for some $c\neq 2$.  Note that in
the case of columns of the form $\stack{n-1}{n}$ (resp.
$\stack{n-1}{\bar{n}}$) the unbracketed $-$ in the $n$-signature
(resp.  $(n-1)$-signature) from $n-1$ cannot be bracketed, since $n$
and $\bar{n}$ may not appear in the same row.

Now, consider a column $\stack{a}{\bar{b}}$ in $T_{2}$, so we have
$a<b$, and thus also $\bar{a}>\bar{b}$.  Recall that $T_{2}$ has no
$n$'s or $\bar{n}$'s, so $b\leq n-1$.  This column contributes $-$'s
to the $a$-signature and the $(b-1)$-signature of $T$.  In this case,
these $-$'s may be bracketed.  Due to the conditions that the rows and
columns of $T$ are increasing, the $-$ from the $a$ can only be
bracketed by an $a+1$ in the bottom row of $T_{1}$ and the $-$ from
the $\bar{b}$ can only be bracketed by a $b$ in the bottom row of
$T_{1}$.  Furthermore, the letter above these must be strictly less
than $a$ and $b-1$, respectively.  By the reasoning in the previous
paragraph, we see that to bracket every $-$ engendered by the column
$\stack{a}{\bar{b}}$ we must have two columns of the form
$\stack{1}{a'}$, with each $a'\neq 2$.

If $k_{3}=1$, $T$ has a column of the form $\stack{a}{\bar{a}}$.  We
have two cases; $2\leq a\leq n-1$, and $a=n$ (resp.  $a=\bar{n}$).  In
the first case, we have a $-$ in the $(a-1)$-signature from the
$\bar{a}$ in this column.  Because of the prohibition against
configurations of the form $\stack{a}{}\stack{a}{\bar{a}}$, this $-$
can only be bracketed by a $+$ from an $a$ in the bottom row of
$T_{1}$.  Therefore, this column engenders another column of the form
$\stack{1}{a'}$.  In the case of $a=n$ (resp.  $a=\bar{n}$), we have a
$-$ in the $(n-1)$-signature (resp.  $n$-signature) which cannot be
bracketed.

To bracket the maximal number of $-$'s (i.e., to minimize $\lan
c,\varphi(T)\ran$) we see that unless $T_{3}=\stack{n}{\bar{n}}$ or
$T_{3}=\stack{\bar{n}}{n}$, we must have
\begin{equation}    \label{eq:t1t2t3}
    T_{1}T_{2}T_{3}=
    \underbrace{
    \stack{1}{b_{1}}
    \cdots
    \stack{1}{b_{2k_{2}+k_{3}}}
    }_{2k_{2}+k_{3}}
    \underbrace{
    \stack{a_{2k_{2}+k_{3}+1}}{b_{2k_{2}+k_{3}+1}}
    \cdots
    \stack{a_{k_{1}}}{b_{k_{1}}}
    }_{k_{1}-(2k_{2}+k_{3})}
    \underbrace{
    \stack{a_{k_{1}+1}}{\bar{b}_{k_{1}+1}}
    \cdots
    \stack{a_{k_{1}+k_{2}+k_{3}}}{\bar{b}_{k_{1}+k_{2}+k_{3}}}
    }_{k_{2}+k_{3}},
\end{equation}
where each column in the first block contributes $3$ to $\lan
c,\varphi(T)\ran$, each column in the second block contributes $2$ to
$\lan c,\varphi(T)\ran$, and the third block contributes nothing.  In
the case $T_{3}=\stack{n}{\bar{n}}$ or $T_{3}=\stack{\bar{n}}{n}$, we
have $k_{2}=0$, so we simply have $T_{1}T_{2}T_{3}=T_{1}T_{3}$, where
each column in $T_{1}$ increases $\lan c,\varphi(T)\ran$ by at least
$2$ and $T_{3}$ increases $\lan c,\varphi(T)\ran$ by $1$.  We
therefore have in the first case $\lan c,\varphi(T)\ran\geq
3(2k_{2}+k_{3})+2(k_{1}-2k_{2}-k_{3})=2k_{1}+2k_{2}+k_{3}$, and in the
second case $\lan c,\varphi(T)\ran\geq
2k_{1}+k_{3}=2k_{1}+2k_{2}+k_{3}$ , as we wished to show.

Since by Lemma \ref{lem:lb} elements in $\tbs$ have level at least
$s$, it follows that when $T\in(\tbs\cap B(s\Lt))_{\min}$, we have
$k_{1}+k_{2}\leq \lfloor s/2\rfloor$, and by $\dual$-duality that
$k_{4}+k_{5}\leq \lfloor s/2\rfloor$.  Furthermore, since
$s=k_{1}+k_{2}+k_{3}+k_{4}+k_{5}$, it follows that
$k_{1}+k_{2}=k_{4}+k_{5}=\lfloor\frac{s}{2}\rfloor$ and $k_{3}=0$ if
$s$ is even and $k_{3}=1$ if $s$ is odd.

\section{Proof of Lemma~\ref{lem:mixed}} \label{app:C}
By using the reverse of Algorithm \ref{alg:iota}, it
suffices to show the following:
    \begin{enumerate}
	\item $T$ has a $1$;
    	\item $T$ has a $\bar{1}$;
    	\item after removing all $1$'s and $\bar{1}$'s, applying the
	Lecouvey $D$ relations will reduce the width of $T$.
    \end{enumerate}
    
The proof of Lemma \ref{lem:plus minus} shows that if $k_{2}\neq 0$,
or $k_{3}=1$ and $T_{3}\neq\stack{n}{\bar{n}}, \stack{\bar{n}}{n}$,
then $T$ has a $1$.  By $\dual$-duality, if $k_{4}\neq 0$, or the same
condition is placed on $k_{3}$ and $T_{3}$, then $T$ has a $\bar{1}$.
We will show that if $k_{2}+k_{3}\neq 0$, then $k_{3}+k_{4}\neq 0$,
which will prove statements (1) and (2) above.

If $k_{3}=1$ this statement is trivial, so we assume $k_{3}=0$.  We
show that the assumptions $k_{2}\neq0$ and $k_{4}=0$ lead to a
contradiction.  From the proof of Lemma \ref{lem:plus minus}, we know
that for $T$ to be minimal, every $-$ from $T_{5}$ must be bracketed.
Because of the increasing conditions on the rows and columns of $T$,
the $-$'s from the bottom row of $T_{5}$ cannot be bracketed by $+$'s
from $T_{5}$, so there must be at least $k_{5}$ $+$'s from
$T_{1}T_{2}$.  Inspection of \eqref{eq:t1t2t3} shows us that the first
block contributes no $+$'s, the second block contributes
$k_{1}-2k_{2}$ many $+$'s, and the third block contributes $2k_{2}$
many $+$'s.  We thus have $k_{1}\geq k_{5}$; but Lemma \ref{lem:plus
minus} tells us that $k_{1}+k_{2}=k_{4}+k_{5}$, contradicting our
assumption that $k_{2}\neq 0$ and $k_{4}=0$.

For the proof of statement (3), we must show that every configuration
$\stack{}{ a}\stack{c}{ b}$ in $T$ avoids the following patterns
(recall the Lecouvey $D$ sliding algorithm from section
\ref{sec:plac}): $\stack{}{x}\stack{n}{\bar{n}} $ and
$\stack{}{x}\stack{\bar{n}}{n} $ with $x\leq n-1$;
$\stack{}{n-1}\stack{\bar{n}}{\overline{n-1}} $ ;
$\stack{}{n-1}\stack{n}{\overline{n-1}} $ ; and $c\geq a$, unless
$c=a=\bar{b}$.  If $T$ has any of these patterns, the top row will not
slide over.

First, simply observe that the first four specified configurations
exclude the possibility of having a column of the form
$\stack{a}{\bar{b}}$ other than $\stack{n}{\bar{n}}$ or
$\stack{\bar{n}}{n}$.  It therefore suffices to show that the presence
of a column $\stack{d}{\bar{e}}$, $2\leq d,e\leq n-1$ implies that $T$
avoids $c\geq a$, unless $c=a=\bar{b}$.  We break our analysis of this
criterion into several special cases:

{\bf Case 1:} $a$ and $b$ are barred, $c$ is unbarred: trivial.

{\bf Case 2:} $a$ is unbarred, $b$ and $c$ are barred: This excludes 
the possibility of having $\stack{d}{\bar{e}}\in T$.

{\bf Case 3:} $a$ and $c$ are unbarred, $b$ is barred: We know the $-$
in the $c$-signature from $c$ must be bracketed; if it is by $b$, we
have $b=\bar{c}$.  As we saw in the proof of Lemma \ref{lem:plus
minus}, we must have the $-$ in the $(c-1)$-signature from $\bar{c}$
bracketed by a $c$ in the bottom row.  This forces $a\geq c$.  If the
$-$ in the $c$-signature from $c$ is bracketed by a $c+1$, it also
must be in the bottom row, forcing $a> c$.

{\bf Case 4:} $a$, $b$, $c$ all unbarred: Suppose $c\geq a$.  Let $d$
be the leftmost unbarred letter weakly to the right of $c$ which does
not have its $-$ bracketed by the letter immediately below it.  (Such
a letter exists, since we assume the occurence of
$\stack{d}{\bar{e}}\in T$, except when $d=e$; this case will be
treated below.)  This letter $d$ must be bracketed by a $d+1$ in the
bottom row, and it must be weakly to the left of $a$.  But we have
$d+1>d\geq c\geq a\geq d+1$; contradiction.

If instead we have a $\stack{d}{ \bar{d}}$ column, we must have the
$-$ from the $\bar{d}$ bracketed by a $d$ in the bottom row weakly to
the left of $a$.  In this case $d\geq c\geq a\geq d$ so $c=d$, and we
have a $\stack{d}{ }\stack{d}{\bar{d}}$ configuration, contradicting
our assumption that $T\in B(s\Lt)$.

{\bf Case 5:} $a$, $b$, $c$ all barred: Similarly to case 4, suppose
$c\geq a$ and let $d$ be the rightmost barred letter weakly to the
left of $a$ which does not have its $+$ bracketed by the letter
immediately above it.  (If none exists, we have a $\stack{\bar{d} }{
d}$ case, see below.)  It must be bracketed by a
$\overline{\bar{d}+1}$ in the top row to the right of $c$.  We then
have $d\leq a\leq c\leq\overline{\bar{d}+1}$; contradiction.

If we have a $\stack{\bar{d}}{d}$ column (note that $d$ is barred),
the $+$ from the $\bar{d}$ must be bracketed by a $d$ in the top row
to the right of $c$.  This implies that $d\geq c\geq a\geq d$, so
$a=d$ and we have a $\stack{\bar{d}}{d}\stack{}{d}$ configuration,
again contradicting our assumption that $T\in B(s\Lt)$.

\section{Proof of Lemma~\ref{lem:unmixed}} \label{app:D}
In the notation of section~\ref{sec:inj}, we have
$k_{2}=k_{4}=0$, and if $k_{3}=1$, $T_{3}=\stack{n}{\bar{n}}$ or
$T_{3}=\stack{\bar{n}}{n}$.  Lemma \ref{lem:plus minus} thus tells us
that $k_{1}=k_{5}$.

Next we show that a column $\stack{j}{i}$ must be of the form
$\stack{i-1}{i}$ for $T$ to be in $\tbs_{\min}$.  For $i=2$ we have
$j=1$ by columnstrictness.  Now suppose that $\stack{j}{i}$ is the
leftmost column such that $j<i-1$.  Then $j$ contributes a $\La_j$ to
$\varphi(T)$ and hence $\lan c,\varphi(T)\ran \ge 2k_{1}+k_{3}+1=s+1$,
so that $T$ is not minimal.  By a similar argument $\lan
c,\varepsilon(T)\ran>s$ unless all columns of the form
$\stack{\bar{i}}{\bar{j}}$ must obey $j=i-1$.

A column $\stack{\overline{i}}{\overline{i-1}}$ for $i>2$ 
(resp. $\stack{n}{\overline{n-1}}$) contributes a $-$ to the 
$(i-2)$-signature (resp. $(n-2)$-signature) of $T$. This $-$ can only be 
compensated by a $+$ in the $(i-2)$-signature (resp. $(n-2)$-signature) from 
a column $\stack{i-1}{i}$ (resp. $\stack{n-1}{\overline{n}}$). 
Hence for $T$ to be minimal the number of columns of the form $\stack{i-1}{i}$ 
(resp. $\stack{n-1}{\overline{n}}$) needs to be the same as the number of columns of 
the form $\stack{\overline{i}}{\overline{i-1}}$ (resp. $\stack{n}{\overline{n-1}}$).
This proves that $T\in \T(s)$.

\end{document}